\documentclass[a4paper,12pt,twoside]{article}
\usepackage{graphicx}\usepackage{a4wide}
\usepackage{amssymb,amsmath,amsthm}
\usepackage{siunitx}
\usepackage{subfigure,color,colordvi,graphicx,xcolor}
\usepackage[only,llbracket,rrbracket,interleave]{stmaryrd}
\usepackage[sort,compress]{cite}
\usepackage{enumitem}
\usepackage{sidecap, caption}

\newcommand{\hl}[1]{{\color{black} #1\color{black}}}

\makeatletter
\newcommand*\wt[2][0.2ex]{%
        \begingroup
        \mathchoice{\wt@helper{#1}{#2}{\displaystyle}{\textfont}}
                   {\wt@helper{#1}{#2}{\textstyle}{\textfont}}
                   {\wt@helper{#1}{#2}{\scriptstyle}{\scriptfont}}
                   {\wt@helper{#1}{#2}{\scriptscriptstyle}{\scriptscriptfont}}%
        \endgroup
        #2%
}
\newcommand*\wt@helper[4]{%
        \def\currentfont{\the#41}%
        \def\currentskewchar{\char\the\skewchar\currentfont}%
        \setbox\tw@\hbox{\currentfont#2\currentskewchar}%
        \dimen@ii\wd\tw@
        \setbox\tw@\hbox{\currentfont#2{}\currentskewchar}%
        \advance\dimen@ii-\wd\tw@
        \rlap{\raisebox{-#1}{$\m@th#3\kern\dimen@ii\widetilde{\phantom{#2}}$}}%
}
\makeatother

\usepackage{booktabs,array,dcolumn}
\newcommand{\PreserveBackslash}[1]{\let\temp=\\#1\let\\=\temp}
\newcolumntype{C}[1]{>{\PreserveBackslash\centering}p{#1}}
\newcolumntype{R}[1]{>{\PreserveBackslash\raggedleft}p{#1}}
\newcolumntype{L}[1]{>{\PreserveBackslash\raggedright}p{#1}}

\newcommand{\bm}[1]{\text{\boldmath $#1$\unboldmath}}

\newcommand{\RR}{\mathbb{R}}
\newcommand{\II}{\mathbb{I}}
\newcommand{\elone}{\ensuremath{\mathcal{L}^1}}

\newcommand{\bx}{\bm{x}}

\newcommand{\bz}{\bm{z}}
\newcommand{\bb}{\bm{b}}
\newcommand{\bd}{\bm{d}}
\newcommand{\bn}{\bm{n}}
\newcommand{\bt}{\bm{t}}
\newcommand{\bhx}{\bm{\widehat{x}}}
\newcommand{\brOne}{\bm{r}_{1,K}}
\newcommand{\brTwo}{\bm{r}_{2,K}}

\newcommand{\OmegaI}{\Omega_{\mathrm{I}}}
\newcommand{\OmegaE}{\Omega_{\mathrm{E}}}
\newcommand{\fI}{f_{\mathrm{I}}}
\newcommand{\fE}{f_{\mathrm{E}}}
\newcommand{\eI}{e_{\mathrm{I}}}
\newcommand{\eE}{e_{\mathrm{E}}}
\newcommand{\chiI}{\chi_{\mathrm{I}}}
\newcommand{\chiE}{\chi_{\mathrm{E}}}
\newcommand{\muI}{\mu_{\mathrm{I}}}
\newcommand{\muE}{\mu_{\mathrm{E}}}
\newcommand{\pI}{\mathrm{p}_{\mathrm{I}}}
\newcommand{\pE}{\mathrm{p}_{\mathrm{E}}}
\newcommand{\ptI}{\mathrm{\tilde{p}}_{\mathrm{I}}}
\newcommand{\ptE}{\mathrm{\tilde{p}}_{\mathrm{E}}}
\newcommand{\PP}{\mathrm{P}}
\newcommand{\lamOne}{\lambda_{1,K}}
\newcommand{\lamTwo}{\lambda_{2,K}}
\newcommand{\dePhi}{\delta\!\phi}

\newcommand{\Heps}{H_{\varepsilon}}
\newcommand{\Ksig}{K_{\sigma}}
\newcommand{\Ktau}{K_{\tau}}

\newcommand{\F}{\mathcal{F}}
\newcommand{\Frsf}{\mathcal{F}_{_{\!\!\mathrm{RSF}}}}
\newcommand{\Fbay}{\mathcal{F}_{_{\!\!\mathrm{BAY}}}}
\newcommand{\G}{\mathcal{G}}
\newcommand{\Grsf}{\mathcal{G}_{_{\!\!\mathrm{RSF}}}}
\newcommand{\Gbay}{\mathcal{G}_{_{\!\!\mathrm{BAY}}}}
\newcommand{\sRSF}{s_{_{\!\mathrm{RSF}}}}
\newcommand{\sBAY}{s_{_{\!\mathrm{BAY}}}}

\newcommand{\Vh}{\mathcal{V}_{\! h}}

\newcommand{\Th}{\mathcal{T}_h}
\newcommand{\Khat}{\widehat{K}}
\newcommand{\DeK}{\Delta_K}
\newcommand{\DeKhat}{\widehat{\Delta}_K}
\newcommand{\BK}{\mat{B}_K}
\newcommand{\LK}{\mat{\Lambda}_K}
\newcommand{\MK}{\mat{M}_K}
\newcommand{\RK}{\mat{R}_K}
\newcommand{\ZK}{\mat{Z}_K}
\newcommand{\wK}{\bm{w}_K}

\newcommand{\PZZ}{\bm{P}_{_{\!\!\!\! \DeK}}}
\newcommand{\GDeK}{\mat{G}_{_{\!\! \DeK}}}
\newcommand{\Metric}{\mat{\mathcal{M}}}
\newcommand{\MetricOld}{\mathcal{M}_{_{\!\mathrm{OLD}}}}
\newcommand{\MetricOmega}{\mathcal{M}_{\!\omega}}
\newcommand{\Jfunc}{\mathcal{J}_K}

\newcommand{\hMin}{\texttt{h}_{\texttt{min}}}
\newcommand{\hMax}{\texttt{h}_{\texttt{max}}}
\newcommand{\nRec}{\texttt{n}_{\texttt{Breg}}}
\newcommand{\nOpt}{\texttt{n}_{\texttt{opt}}}
\newcommand{\nAdapt}{\texttt{n}_{\texttt{adapt}}}
\newcommand{\tOpt}{\ensuremath{\overline{\texttt{T}}_{\texttt{opt}}}}
\newcommand{\tAdapt}{\ensuremath{\overline{\texttt{T}}_{\texttt{adapt}}}}

\newcommand{\mat}[1]{\mathbf{#1}}

\newcommand{\grad}{\bm{\nabla}}

\newcommand{\eltwo}{\ensuremath{\mathcal{L}_2}}

\newcommand{\numel}{\ensuremath{\texttt{n}_{\texttt{el}}}}

\newenvironment{keywords}{\begin{quote}\emph{\textbf{Keywords:}}}{\end{quote}}
\newtheorem{remark}{Remark}
\newtheorem{proposition}{Proposition}

\RequirePackage{algorithm}[0.1] 
\usepackage{algorithmic} 

\graphicspath{{figsPDF/}{figsPNG/}{figs/}}

\usepackage{hyphenat}

\begin{document}
\title{Anisotropic mesh adaptation for region-based segmentation accounting for image spatial information}

\author{
\renewcommand{\thefootnote}{\arabic{footnote}}
			  Matteo Giacomini\footnotemark[1]\textsuperscript{ \ ,}\footnotemark[2]\textsuperscript{ \ ,}*  \  and
			  Simona Perotto\footnotemark[3]
}

\date{\today}
\maketitle

\renewcommand{\thefootnote}{\arabic{footnote}}

\footnotetext[1]{Laboratori de C\`alcul Num\`eric (LaC\`aN), ETS de Ingenieros de Caminos, Canales y Puertos, Universitat Polit\`ecnica de Catalunya, Barcelona, Spain.}
\footnotetext[2]{Centre Internacional de M\`etodes Num\`erics en Enginyeria (CIMNE), Barcelona, Spain.}
\footnotetext[3]{MOX, Dipartimento di Matematica, Politecnico di Milano, Piazza L. da Vinci 32, I-20133 Milano, Italy.
\vspace{5pt}\\
* Corresponding author: Matteo Giacomini. \textit{E-mail:} \texttt{matteo.giacomini@upc.edu}
}

\begin{abstract}
A finite element-based image segmentation strategy enhanced by an anisotropic mesh adaptation procedure is presented. The methodology relies on a split Bregman algorithm for the minimisation of a region-based energy functional and on an anisotropic recovery-based error estimate to drive mesh adaptation.
More precisely, a Bayesian energy functional is considered to account for image spatial information, ensuring that the methodology is able to identify inhomogeneous spatial patterns in complex images.
In addition, the anisotropic mesh adaptation guarantees a sharp detection of the interface between background and foreground of the image, with a reduced number of degrees of freedom.
The resulting split-adapt Bregman algorithm is tested on a set of real images showing the accuracy and robustness of the method, even in the presence of Gaussian, salt and pepper and speckle noise.
\end{abstract}


\newpage

\section{Introduction}
\label{sc:Intro}

Image segmentation aims to identify a partition of a given image into subregions characterised by different pixel intensities or textures. This problem is of great interest for several applications, encompassing medical imaging~\cite{Brox-RFB-15}, autonomous vehicles~\cite{Ha-HWKUH-17}, agriculture~\cite{Sammouda-SATA-14} and forensics~\cite{Bazen-BG-01}.

Many techniques have been proposed in the literature to solve image segmentation problems, including thresholding~\cite{Otsu-79}, edge detection~\cite{Ziou-ZT-98}, region-growing~\cite{Zhao-ZWWS-14} and neural networks~\cite{Milletari-MNA-16,Dolz-DDA-18} approaches. In the seminal work~\cite{MumfordShah-89}, D.B. Mumford and J. Shah proposed a variational framework for image segmentation introducing an energy functional whose minimum is associated with the boundary of the object to be segmented. Stemming  from this result, different strategies were  proposed to approximate the Mumford-Shah functional. Major contributions encompass the Chan-Vese approach~\cite{ChanVese-01} based on an active contour model without edge detection, the Ambrosio-Tortorelli functional~\cite{AmbrosioTortorelli-90} and other discrete models~\cite{Chambolle-95,Chambolle-CD-99,Bourdin-99} providing sequences of approximations $\Gamma$-convergent to the Mumford-Shah functional.


The present work focuses on the finite element approximation of region-based image segmentation strategies originated from the Chan-Vese model. These methods rely on the solution of a partial differential equation to describe the evolution of a level-set function associated with the contour of the region to be segmented~\cite{Vese-VC-02,Li-LXGF-05,Wang-WHX-10}. In this context, the treatment of image inhomogeneities and noise represents a crucial aspect for the robustness of the segmentation algorithm. In~\cite{Li-LKGD-08}, the region-scalable fitting energy (RSFE) model was proposed to deal with image inhomogeneities by introducing appropriate kernel functions in the data fitting term, thus accounting for the information provided by the pixel intensity in local regions at a controllable scale. To fully exploit the richness in spatial information of the image, see e.g.~\cite{Zhang-ZLCCL-22}, statistical models leveraging nonlocal information via appropriate estimates of the probability distributions of the pixel intensities were proposed in~\cite{Fisher-KFYCW-05,Brox-BC-09,Zhang-ZZLZ-15}. More recently, bias correction models~\cite{Dong-DJW-19,Weng-WDL-21} and super-resolution techniques~\cite{Calatroni-CCP-21} were introduced to ensure a robust treatment of inhomogeneous images with the goal of  achieving high-resolution segmentations, also in the presence of noise. A known issue, common to the above mentioned models, is the ill-posedness of the associated minimisation problem. To remedy this issue, a strategy to convexify the energy functional was proposed in~\cite{Chan-CEN-06}. In addition, the split Bregman algorithm~\cite{Osher-GO-09} was specifically developed to tackle optimisation problems with regularised functionals as the ones arising in image segmentation~\cite{Osher-YLKO-10,Yang-YZW-13}.

When dealing with the above mentioned numerical optimisation problems, the assessment of the quality of the computed solution (or of quantities of interest depending on the solution) represents a crucial aspect to provide verified and reliable results. In the context of finite element approximations, \emph{a posteriori} error estimates have been successfully developed to provide either quantitative error bounds~\cite{Destuynder-DM-99,Pares-PDH-06,Ern-EV-15,Prudhomme-OP-01,Rodenas-DRZ-07,Pares-PDH-09,Ainsworth-AR-12,Prudhomme-MP-15,Giacomini-GPT-17,Giacomini-18} or qualitative information to drive mesh adaptation strategies~\cite{Babuska-BR-78,Rannacher-BR-13,Huerta-HRDS-99}. In the latter case, particular attention has been devoted to the local control of the geometric features of the mesh (i.e., shape, size and orientation of the elements), by means of anisotropic elements stretched along the direction orthogonal to the gradient of the quantity of interest~\cite{Stein-SO-99,Fortin-DVBFH-02,Perotto-MP-06,Alauzet-LDA-10,Perotto-FMP-11,Perotto-MP-11,Perotto-PPB-12,Fortin-BFB-14,Perotto-AFMP-15,Shontz-KMS-15,Ferro-FMP-18,Perotto-MPS-19}. It is worth noticing that the coupling of finite element approximations and mesh adaptation in the context of image segmentation was previously explored in~\cite{Yaacobson-YG-98,Chambolle-BC-00,Nochetto-DMN-08} and~\cite{Coupez-ZCDJCS-16,Perotto-CMPP-19,Perotto-CFMMNP-20} in isotropic and anisotropic settings, respectively. Nonetheless, all the above models were devised for images featuring sharp interfaces and regions with homogeneous intensities.

In order to treat image inhomogeneities, the present work couples a region-based Bayesian segmentation model, able to account for spatial information, with an anisotropic mesh adaptation procedure. More precisely, the split Bregman method for the segmentation is enriched with a recovery-based error estimate to construct anisotropic meshes guaranteeing an accurate description of the contour of the image, with a reduced number of degrees of freedom. The resulting algorithm is thus capable of exploiting the spatial information to identify complex inhomogeneous patterns in real images, even in the presence of different sources of noise, such as Gaussian, salt and pepper and speckle noise. The remainder of the article is organised as  follows. Section~\ref{sc:ImgSeg} introduces the split Bregman method for a region-based segmentation model based on a Bayesian approach and able to account for the spatial information in the image. In section~\ref{sc:Adaptivity}, the anisotropic recovery-based error estimate is presented along with the split-adapt Bregman method resulting from its coupling with the selected optimisation algorithm. Numerical validation of the split-adapt Bregman algorithm and its application to the segmentation of real images corrupted by different types of noise are presented in section~\ref{sc:Simulations}. Finally, section~\ref{sc:Conclusion} summarises the presented results, \ref{sc:appRSFE} reports the extension of the discussed framework to the RSFE model and \ref{sc:appParam} provides some details about the parameters employed in the numerical simulations.

\section{The optimisation problem of region-based image segmentation}
\label{sc:ImgSeg}

In this section, the formulation of a region-based segmentation model accounting for image spatial information is introduced and the split Bregman algorithm~\cite{Osher-GO-09} to solve the resulting minimisation problem is recalled.

Consider an image $\Omega {\subset} \RR^2$ and denote by $U {:} \Omega {\rightarrow} [0,255]$ the grey level intensity of the associated pixels. A level-set function $\phi$ is employed to describe the unknown contour separating $\Omega$ into two disjoint regions $\OmegaI$ and $\OmegaE$ defined as
\begin{subequations}
\begin{gather}
\label{eq:Regions}
\OmegaI := \{ \bx \in \Omega \ : \ \Heps(\phi(\bx)) = 1 \}  ,\\
\OmegaE := \{ \bx \in \Omega \ : \ \Heps(\phi(\bx)) = 0 \}  ,
\end{gather}
\end{subequations}
where $\Heps$ denotes the regularised Heaviside function
\begin{equation}
\label{eq:Heaviside}
\Heps(\phi) := \frac{1}{2}  + \frac{1}{\pi} \arctan{\left(\frac{\phi}{\varepsilon}\right)} , \quad \varepsilon > 0 .
\end{equation}

In this framework, the image segmentation problem can be written as an optimisation procedure consisting of minimising a suitable energy functional $\F(\phi(\bx))$. In order to enforce the existence of a global minimum for the image segmentation functional $\F$, a constraint on the level-set function is commonly introduced~\cite{Osher-YLKO-10}, namely
\begin{equation}
\label{eq:ConstraintPhi}
-\alpha \leq \phi(\bx) \leq \alpha , \ \text{ for a.e. } \bx \in \Omega .
\end{equation}
In the present work, the level-set function is described by means of a piecewise constant approximation such that
\begin{equation}
\label{eq:Phi}
\phi(\bx) = \begin{cases}
\alpha, & \bx \in \OmegaI \\
-\alpha, & \bx \in \OmegaE \\
\end{cases}  
\end{equation}
and the value $\alpha {=} 1$ is selected for the numerical simulations in section~\ref{sc:Simulations}.

\subsection{A Bayesian approach accounting for spatial information}
\label{sc:Bayes}

In order to handle inhomogeneous images by exploiting the intrinsic spatial information, a region-based Bayesian segmentation model is considered. This model is obtained by maximising the \emph{a posteriori} probability $\PP(\bx {\in} \Omega_i \ | \ U(\bx) {=} \kappa)$ of correctly assigning a given pixel $\bx$ to the region $\Omega_i$ it belongs to, knowing the value $\kappa$ of the pixel intensity $U(\bx)$. From Bayes' theorem, it follows
\begin{equation}
\label{eq:MAP}
\PP(\bx \in \Omega_i \ | \ U(\bx) = \kappa) \propto \mathrm{p}_i(\kappa) \PP(\bx \in \Omega_i) ,
\end{equation}
where the likelihood $\mathrm{p}_i(\kappa) := \PP(U(\bx) {=} \kappa \ | \ \bx {\in} \Omega_i )$ represents the probability that the pixel intensity at $\bx$ is equal to $\kappa$, knowing that point $\bx$ belongs to the region $\Omega_i$, whereas $\PP(\bx {\in} \Omega_i)$ is the prior probability distribution encapsulating any \emph{a priori} information available on the image.

To determine the likelihood, the spatial information of the image is exploited by performing a kernel density estimate~\cite{Parzen-62} and a probability density function (PDF) of the pixel intensity in each region is constructed~\cite{Fisher-KFYCW-05,Brox-BC-09}, namely
\begin{equation}
\label{eq:pIntExt}
\mathrm{p}_i(\kappa) :=  \int_0^{255}{\Ktau(\kappa - \xi) \II_i(\xi) d\xi} ,
\end{equation}
where $\Ktau$ is a Gaussian kernel with standard deviation $\tau$, namely
\begin{equation}
\label{eq:GaussKernel}
\Ktau(\bx) := \frac{1}{2 \pi \tau^2} e^{- \frac{\| \bx \|^2}{2 \tau^2}} ,
\end{equation}
$\| \bx \|$ being the Euclidean norm of the position vector, whereas $\II_i$ is the discrete histogram accounting for the frequency of recurrence of the pixel intensities in the region $\Omega_i$. More precisely, the value $\II_i(\kappa)$ of the recurrence of the intensity $\kappa$ in the region $\Omega_i$ is obtained by computing the cardinality of the set $\{ \bx {\in} \Omega_i  :  U(\bx) {=}\kappa \}$. The procedure in equation~\eqref{eq:pIntExt} thus performs a smoothing of the resulting discrete histogram via a convolution product with the Gaussian kernel $\Ktau$ in order to obtain a continuous definition of the probability distribution $\mathrm{p}_i$.

As mentioned above, the prior $\PP(\bx {\in} \Omega_i)$ can be employed to incorporate any available information on the image into the model, e.g., by exploiting the knowledge on the nature of the noise due to the measurement devices used for image acquisition~\cite{Burger-STJB-13}. In this work, a non-informative shape prior accounting only for the information on the regularity and the length of the contour is considered, namely
\begin{equation}
\label{eq:prior}
\PP(\bx \in \Omega_i) \propto e^{-\nu \| \phi(\bx) \|_{\textrm{TV},g}} ,
\end{equation}
where $\nu$ is a positive regularisation parameter and $\| \phi(\bx) \|_{\textrm{TV},g}$ denotes the weighted total variation (TV) norm, i.e.,  
\begin{equation}
\label{eq:TVnorm}
\| \phi(\bx) \|_{\textrm{TV},g} := \int_{\Omega}{g(U(\bx)) | \nabla \phi(\bx) | d\bx} ,
\end{equation}
defined as the $\elone(\Omega)$ norm of the gradient of the level-set function weighted by the so-called \emph{edge detector} function~\cite{Li-LKGD-08,Osher-YLKO-10}, namely
\begin{equation}
\label{eq:edgeDet}
g(U) = \frac{1}{1 + \beta \| \nabla U \|_2^2} , \quad \beta > 0 ,
\end{equation}
where $\| \cdot \|_2$ denotes the classical $\eltwo(\Omega)$ norm and $\beta$ is a positive scaling factor.

Assuming independence of the pixel intensities at different points $\bx$, see e.g.~\cite{Brox-BC-09}, the \emph{a posteriori} probability~\eqref{eq:MAP} can be rewritten as
\begin{equation}
\label{eq:MAP-prod}
\PP(\bx \in \Omega_i \ | \ U(\bx) = \kappa) \propto \prod_{\bx \in \OmegaI} \pI(U(\bx)) \prod_{\bx \in \OmegaE} \pE(U(\bx)) \ e^{-\nu \| \phi(\bx) \|_{\textrm{TV},g}} .
\end{equation}
It is worth recalling that the maximisation of the functional in equation~\eqref{eq:MAP-prod} is equivalent to the minimisation of its negative logarithm~\cite{Zhu-ZY-96}. Hence, the energy functional for the Bayesian model is given by 
%
%
\begin{equation}
\label{eq:funcBayes}
\begin{aligned}
\Fbay(\phi,\ptI,\ptE)  := &-\int_{\Omega}{ [ \log(\ptI(U(\bx)) + \log(\ptE(U(\bx)) ] d\bx} \\
&+ \nu \int_{\Omega}{g(U(\bx)) | \nabla \phi(\bx) | d\bx} .
\end{aligned}
\end{equation}
The first term drives the separation of $\Omega$  into two regions, $\OmegaI$ and $\OmegaE$, by exploiting the image spatial information encapsulated in the PDF of the pixel intensity
\begin{equation}
\label{eq:probReg}
\mathrm{\tilde{p}}_i(U(\bx)) := \begin{cases}
\mathrm{p}_i(U(\bx)) , & \text{for $\bx \in \Omega_i$} \\
\zeta , & \text{elsewhere}
\end{cases}
\end{equation}
where $\zeta {>} 0$ is a regularisation parameter which guarantees that $\log(\mathrm{\tilde{p}}_i(U(\bx)))$ is well-defined. 
The second term in equation~\eqref{eq:funcBayes} features the weighted TV norm of the level-set function and it is responsible for the regularisation of the segmented contour by enhancing its smoothness~\cite{Osher-BEVTO-07}. It is worth noticing that the introduction of this regularisation term in equation~\eqref{eq:funcBayes} has two goals. On the one hand, it fosters the detection of local changes in the pixel intensity through the edge detection function defined in~\eqref{eq:edgeDet}. On the other hand, it contributes to the convexification of the energy functional, thus ensuring the existence of a global minimum~\cite{Chan-CEN-06, Osher-YLKO-10}.

\begin{remark}[Regularisation term]
It is well-known that the TV norm $\| \phi(\bx) \|_{\textrm{TV},g}$ denotes a total variation measure. If $\Omega$ is a Caccioppoli set, such a measure represents the perimeter of the set $\Omega$ itself~\cite{Giusti-GW-84}. In this context, the regularisation term in equation~\eqref{eq:funcBayes} can be interpreted as a penalisation on the length of the boundary, preventing high-frequency oscillations of the segmented contour.
\end{remark}

\subsection{Split Bregman algorithm for image segmentation}
\label{sc:SplitBreg}

The minimisation of the Bayesian energy functional introduced in equation~\eqref{eq:funcBayes} can be recast as
\begin{equation}
\label{eq:minL1phi}
\min_{\phi} \Big\{ \G(\phi(\bx)) + \nu \| \phi(\bx) \|_{\textrm{TV},g} \Big\} ,
\end{equation}
where the functional $\G(\phi(\bx))$ is given by
\begin{equation}
\label{eq:funcG}
\G(\phi(\bx)  = \Gbay(\phi(\bx)) := -\int_{\Omega}{ [ \log(\ptI(U(\bx)) + \log(\ptE(U(\bx)) ] d\bx} .
\end{equation}

It is worth recalling that the weighted TV norm of the level-set function in equation~\eqref{eq:minL1phi} can be equivalently interpreted as a weighted $\elone(\Omega)$ norm of its gradient, that is $\| \grad \phi \|_{1,g} := \| \phi \|_{\textrm{TV},g}$. Moreover, following~\cite{Zhang-WYYZ-08}, an auxiliary variable $\bd {=} \grad \phi$ is introduced and problem~\eqref{eq:minL1phi} is rewritten as
\begin{equation}
\label{eq:minL1phiD}
\min_{\phi, \bd} \Big\{ \G(\phi(\bx)) + \nu \| \bd(\bx) \|_{1,g} +  \frac{\mu}{2} \| \bd(\bx) - \grad \phi(\bx) \|_2^2 \Big\} ,
\end{equation}
where the last term enforces the constraint $\bd  {=} \grad \phi$ with the penalty parameter $\mu {>} 0$.

\hl{
The main difficulty of the minimisation problem~\eqref{eq:minL1phiD} is represented by the $\elone(\Omega)$ regularisation term in the objective functional.  In this work, the split Bregman method is employed to solve problem~\eqref{eq:minL1phiD}, according to a \emph{staggered} rationale~\cite{Osher-GO-09}. In particular, the Bregman splitting allows to account for the $\elone(\Omega)$ and the $\eltwo(\Omega)$ components in the energy functional independently, thus providing an iterative strategy where minimisation with respect to $\phi$ and $\bd$ is performed separately. Besides reducing the number of globally coupled unknowns,  the advantage of this approach stems from the simple structure of the two resulting subproblems: on the one hand,  $\phi$ is not affected by the $\elone(\Omega)$ component of the functional, so that a derivative-based (here, a gradient flow) strategy can be employed for its computation; on the other hand,  a \emph{shrinkage} operator is suitable to efficiently determine $\bd$. Of course, alternative techniques not relying on the Bregman splitting could be considered to solve problem~\eqref{eq:minL1phiD}, such as variational methods exploiting an alternating minimisation algorithm~\cite{Bourdin-99, Perotto-CMPP-19}. However, this issue lies beyond the scope of the present work.
}

At iteration $k$, the split Bregman method, detailed in algorithm~\ref{alg:splitBregman}, performs the following three steps:
\begin{subequations}
\label{eq:splitBregman}
\begin{align}
\text{\underline{Step A:} } \phi^{k+1} &= \min_{\phi} \Big\{ \G(\phi) +  \frac{\mu}{2} \| \bd^k - \grad \phi - \bb^k \|_2^2 \Big\}  ,
\label{eq:splitBregman-A}\\
\text{\underline{Step B:} }\bd^{k+1} &= \min_{\bd} \Big\{ \nu \| \bd \|_{1,g} +  \frac{\mu}{2} \| \bd - \grad \phi^{k+1} - \bb^k \|_2^2 \Big\} ,
\label{eq:splitBregman-B} \\
\text{\underline{Step C:} } \bb^{k+1} &= \bb^k + \grad \phi^{k+1} - \bd^{k+1} .
\label{eq:splitBregman-C}
\end{align}
\end{subequations}

\begin{algorithm}[!ht]
\caption{Split Bregman algorithm for image segmentation}\label{alg:splitBregman}
\begin{algorithmic}[1]
\REQUIRE{Tolerance $\eta^\star$ for the stopping criterion and initial contour $\phi_0$.}
\STATE{Construct  the computational mesh $\Th$.}
\STATE{Set $k {=} 0$, $\phi_h^0 {=} \phi_0$, $\bd_h^0 {=} \grad \phi_0$ and $\bb_h^0 {=} \bm{0}$.}
\WHILE{$\left\| \phi_h^{k+1} - \phi_h^k \right\|_2 > \eta^\star \left\| \phi_h^k \right\|_2$}
\STATE{Compute the contour $\phi_h^{k+1}$ solving equation~\eqref{eq:evolutionPhiFE}.}
\STATE{Compute the auxiliary variable $\bd_h^{k+1}$ via the shrinkage operator~\eqref{eq:shrinkD}.}
\STATE{Compute the Bregman update $\bb_h^{k+1}$ using equation~\eqref{eq:splitBregman-C}.}
\STATE{$k \gets k+1$.}
\ENDWHILE
\ENSURE{Boundary $\phi_h^k$ of the segmented region.}
\end{algorithmic}
\end{algorithm}

Starting from a guess of the contour, $\phi {=} \phi_0$, the variables $\bd^0 {=} \grad \phi_0$ and $\bb^0 {=} \bm{0}$ are set.  Step A is responsible for the computation of $\phi^{k+1}$ after fixing the last available approximations for $\bd$ and $\bb$, obtained from iteration $k$. Following the gradient flow equation approach proposed in~\cite{Chan-CEN-06}, the solution $\phi^{k+1}$ of the minimisation problem~\eqref{eq:splitBregman-A} is computed as the steady-state solution of the evolution equation
\begin{equation}
\label{eq:evolutionPhi}
\left\{
\begin{aligned}
\frac{\partial \phi}{\partial t} - \Delta \phi &= -\frac{1}{\mu} s^k + \grad \cdot (\bb^k - \bd^k) , && \text{ in $\Omega$, $t>0$} \\
\grad \phi \cdot  \bn &= 0 && \text{ in $\partial\Omega$, $t>0$} \\
\phi &= \phi_0, && \text{ in $\overline{\Omega}$, $t=0$} 
\end{aligned}
\right.
\end{equation}
where the source term $s^k$ is obtained by evaluating the expression
\begin{equation}
\label{eq:sourceTerm}
s(\bx) = \sBAY(\bx) := 
-\log(\ptI(U(\bx)) - \log(\ptE(U(\bx)) 
\end{equation}
using the last computed level-set function $\phi^k$.

\hl{
In this work, a low-order conforming finite element framework is considered for the spatial approximation of the level-set function representing the boundary of the region to be segmented. In this context,  the accuracy of the segmented boundary strongly depends on the resolution of the underlying computational mesh.  To improve the quality of such a representation, while limiting the additional computational cost of the procedure, adapted meshes are employed. Of course, the accuracy of the solution could be improved by means of high-order discretisation schemes and high-order meshes, as recently explored by several authors (see, e.g.,~\cite{Giani-AGC-19,Ramanuj-RS-19,Falcone-FPT-20}).  Although high-order schemes allow to achieve accurate level-set approximations albeit using coarse meshes,  it is well-known that such methods are less robust than low-order approaches, especially when pure convection phenomena are considered. More precisely,  high-order methods might require the introduction of artificial viscosity or slope limiters~\cite{CockburnShu-89,PerssonPeraire-06,Casoni-HCP-12} to avoid non-physical oscillations in the presence of steep gradients, as the ones in the neighbourhood of the boundary to be segmented. These considerations justify the proposed choice of a low-order, nevertheless robust, spatial solver.}

\hl{To obtain the fully-discrete form of equation~\eqref{eq:evolutionPhi},  an appropriate time integration scheme needs to be introduced. As previously mentioned, the solution of problem~\eqref{eq:splitBregman-A} is computed as the steady-state solution of the evolution equation~\eqref{eq:evolutionPhi}. In this context, transient effects are negligible and the time marching scheme is employed as a relaxation approach to speed-up the convergence of the iterative algorithm. Given an artificial time step $\Delta t$, the time derivative can be discretised using the backward Euler scheme, that is,
$$
\frac{\partial}{\partial t}\phi(t^{k+1}) \simeq \frac{1}{\Delta t} (\phi_h^{k+1} -\phi_h^k) .
$$
Other numerical approximations of equation~\eqref{eq:evolutionPhi} can be derived using high-order time integrators such as backward difference formulae or Runge-Kutta methods~\cite{Hairer-book}. Nonetheless,  the additional accuracy provided by these schemes is mainly of interest in the simulation of transient phenomena. This justifies the employment of the backward Euler scheme to approximate the time dependence in equation~\eqref{eq:evolutionPhi}. Moreover,  in the numerical assessment of section~\ref{sc:Simulations}, it was verified that the backward Euler method converges in few iterations (in general, at most five) at each iteration of the split Bregman algorithm.}

\hl{To construct the discretisation of step A of the split Bregman algorithm,  first, a computational mesh $\Th$ consisting of $\numel$ conforming simplicial elements, $K_i$, $i {=} 1,\ldots,\numel$, is defined (Algorithm~\ref{alg:splitBregman} - Step 1). }Let $\Vh$ be the space of piecewise linear continuous finite element functions on the mesh $\Th$ and denote by $\phi_h$, $\bd_h$ and $\bb_h$ the finite element approximations of $\phi$, $\bd$ and $\bb$, respectively. The discrete counterpart of equation~\eqref{eq:evolutionPhi}, obtained by considering a conforming finite element discretisation in space and an implicit Euler scheme for time marching, is: find $\phi^{k+1} {\in} \Vh$ such that
\hl{
\begin{equation}
\label{eq:evolutionPhiFE}
\begin{aligned}
\int_{\Omega}{\dePhi \, \phi_h^{k+1} d\bx} &+ \Delta t \int_{\Omega}{\grad \dePhi \cdot \grad \phi_h^{k+1} d\bx} \\
&= \int_{\Omega}{\dePhi \, \phi_h^k d\bx} + \Delta t \int_{\Omega}{\dePhi \left(-\frac{1}{\mu} s^k + \grad \cdot (\bb_h^k - \bd_h^k) \right) d\bx} ,
\end{aligned}
\end{equation}
}
for all $\dePhi {\in} \Vh$ (Algorithm~\ref{alg:splitBregman} - Step 4).

The value $\phi_h^{k+1}$ yielded by the discrete problem~\eqref{eq:evolutionPhiFE} is thus used in step B to determine a new approximation for $\bd_h$ via the solution of the minimisation problem~\eqref{eq:splitBregman-B} (Algorithm~\ref{alg:splitBregman} - Step 5). This is achieved via a thresholding operation, namely
\begin{equation}
\label{eq:shrinkD}
\bd_h^{k+1} = \texttt{shrink} \left\{ \bb_h^k + \grad \phi_h^{k+1} , \frac{\nu}{\mu} g(U) \right\} ,
\end{equation}
where $\texttt{shrink}$ denotes the \emph{shrinkage} operator~\cite{Osher-GO-09} defined as
\begin{equation}
\label{eq:shrinkOp}
\texttt{shrink} \left\{ f , \gamma \right\} := \frac{f}{|f|} \max \{ |f| - \gamma , 0 \} .
\end{equation}

Finally, step C of the split Bregman method computes the correction $\bb_h^{k+1}$ accounting for the discrepancy between $\grad \phi_h^{k+1}$ and $\bd_h^{k+1}$ by means of equation~\eqref{eq:splitBregman-C} (Algorithm~\ref{alg:splitBregman} - Step 6), before the counter of the iterative procedure is updated (Algorithm~\ref{alg:splitBregman} - Step 7) and a new  iteration starts. 

The split Bregman method stops when a steady-state solution is achieved, that is, when the relative variation of the level-set function from iteration $k$ to iteration $k {+} 1$ is below a user-defined tolerance $\eta^\star$ (Algorithm~\ref{alg:splitBregman} - Step 3). Of course, other stopping criteria may also be devised to exploit specific information provided by the model. For instance, a spatial function measuring the discrepancy between $\ptI(U(\bx))$ and $\ptE(U(\bx))$ at each point $\bx$ could be considered, namely
\begin{equation}
\label{eq:probVariation}
\delta_{\mathrm{p}}(\bx) := \left| \ptI(U(\bx)) - \ptE(U(\bx)) \right| .
\end{equation}
The resulting stopping criterion accounting for the variation of the probability distributions from iteration $k$ to iteration $k {+} 1$ would thus follow by substituting step 3 of algorithm~\ref{alg:splitBregman} with the verification of the condition $\| \delta_{\mathrm{p}}^{k+1} - \delta_{\mathrm{p}}^k  \|_2 > \eta^\star \| \delta_{\mathrm{p}}^k \|_2$.

It is worth noticing that the described split Bregman strategy is suitable also for other image segmentation models that can be formulated according to equation~\eqref{eq:minL1phi}. This is the case of the RSFE functional~\cite{Li-LKGD-08,Osher-YLKO-10} employed for benchmarking in some simulations of section~\ref{sc:Simulations}. A presentation of the RSFE model is reported in~\ref{sc:appRSFE}, whereas the details of the parameters involved in the definition of  the associated energy functional and in the construction of the split Bregman algorithm are provided in~\ref{sc:appParam}.

\section{Anisotropic mesh adaptation enhancing image segmentation}
\label{sc:Adaptivity}

In this section, a mesh adaptation strategy is introduced into the split Bregman algorithm. Since a piecewise constant definition is employed for the level-set function $\phi$, see equation~\eqref{eq:Phi}, the region where $\grad \phi$ is different from zero is associated with the neighbourhood of the contour to be segmented. Hence, following the seminal works~\cite{ZZ-1987,ZZ-92-part1,ZZ-92-part2} by O.C. Zienkiewicz and J.Z. Zhu on recovery-based error estimates, a local error indicator is constructed to evaluate the error of the gradient of the level-set function measured in the $\eltwo(\Omega)$ norm. The computation is carried out in an anisotropic framework~\cite{Perotto-FP-01,Perotto-FP-03} in order to exploit the intrinsic directionality of the gradient of the level-set function. The resulting strategy allows to increase the accuracy of the segmented region while reducing the number of degrees of freedom required for the description of the contour, as shown by the numerical experiments in section~\ref{sc:Simulations}.

\subsection{The anisotropic setting}
\label{sc:aniso}

According to the anisotropic setting in~\cite{Perotto-FP-01,Perotto-FP-03}, the geometric properties (i.e., shape, size and orientation) of each mesh element $K {\in} \Th$ are described through the spectral properties of the affine invertible map $T_K {:} \Khat {\rightarrow} K$. This map transforms the reference equilateral triangle $\Khat$ inscribed in the unit circle with centre at the origin into a generic triangle $K$ inscribed in an ellipse. Any point $\bx {=} (x,y)^T {\in} K$ can thus be expressed as a function of the coordinates $\bhx {=} (\widehat{x},\widehat{y})^T$ in the reference element $\Khat$, being
\begin{equation}
\bx = T_K(\bhx) := \MK \bhx +\wK ,
\label{eq:Map}
\end{equation}
where $\MK {\in} \RR^{2 \times 2}$ is the Jacobian of the transformation responsible for the rotation and deformation of the reference element $\Khat$, whereas $\wK {\in} \RR^2$ is a shift vector accounting for rigid translations.

The polar decomposition $\MK {=} \BK \ZK$ of the Jacobian allows to identify a symmetric positive definite matrix $\BK {\in} \RR^{2 \times 2}$ describing the deformation of the triangle $K$ and an orthogonal matrix $\ZK {\in} \RR^{2 \times 2}$ accounting for the rotation of $K$. Moreover, from the spectral decomposition of $\BK$, it follows that $\BK {=} \RK^T \LK \RK$, where $\RK^T := [\brOne,\brTwo]$ is the matrix of the right eigenvectors and $\LK := \operatorname{diag}(\lamOne,\lamTwo)$ is the diagonal matrix of the corresponding eigenvalues, with $\lamOne \geq \lamTwo > 0$. The eigenvectors $\brOne$ and $\brTwo$ identify the directions of the semi-axes of the ellipse circumscribed to the triangle $K$, whereas $\lamOne$ and $\lamTwo$ measure the corresponding lengths. Finally, the quantity $s_K := \lamOne / \lamTwo \geq 1$, known as aspect ratio or stretching factor, quantifies the elemental anisotropy of the triangle $K$. In particular, $s_K  {=} 1$ corresponds to the equilateral shape (i.e., the isotropic case), whereas the higher the  value of $s_K$, the more stretched the resulting triangle.

\subsection{An anisotropic recovery-based error estimate}
\label{sc:ZZ}

In order to define a recovery-based error estimate, first the  so-called recovered gradient needs to be constructed. Following~\cite{Perotto-MP-10, Perotto-PPB-12}, a piecewise constant reconstruction of the discrete gradient $\grad \phi_h$ over the patch $\DeK := \{ T {\in} \Th \ : \ T  {\cap} K {\neq} \emptyset \}$, associated with element $K$, is considered. More precisely, the average of the gradient of the discrete solution $\phi_h$ over the patch is computed, weighted by the area of the patch elements. The resulting recovered gradient is
\begin{equation}
\label{eq:recoveryConst}
\PZZ \phi_h (\bx) := \frac{1}{| \DeK |} \sum_{T \in \DeK} | T | \left. \grad \phi_h \right|_T , \quad \text{for } \bx \in K .
\end{equation}

Thus, the anisotropic recovery-based error estimate~\cite{Perotto-MP-10} is given by
\begin{equation}
\label{eq:estimate}
\eta^2 = \sum_{K \in \Th} \eta_K^2 ,
\end{equation}
where $\eta_K$ denotes the contribution of element $K$  to the global estimate, with
\begin{equation}
\label{eq:estimateLoc}
\eta_K^2 := \frac{1}{\lamOne \lamTwo} \left[ \lamOne^2 \brOne^T \GDeK \brOne + \lamTwo^2 \brTwo^T \GDeK \brTwo \right] ,
\end{equation}
$\GDeK {\in} \RR^{2 \times 2}$ being the symmetric positive semidefinite matrix with components
\begin{equation}
\label{eq:GDeltaK}
\left[ \GDeK \right]_{ij} := \sum_{T \in \DeK} \int_T{ \left[\PZZ \phi_h -  \left. \grad \phi_h \right|_{\DeK} \right]_i \left[\PZZ \phi_h -  \left. \grad \phi_h \right|_{\DeK} \right]_j d\bx} .
\end{equation}
It is worth noticing that the scaling factor $(\lamOne \lamTwo)^{-1}$ in equation~\eqref{eq:estimateLoc} guarantees the consistency of the anisotropic error estimate with the isotropic case characterised by $\lamOne {=} \lamTwo$.

Finally, the estimate in equation~\eqref{eq:estimateLoc} can be rewritten by means of a geometric scaling as
\begin{equation}
\label{eq:estimateLocRescale}
\eta_K^2 := \lamOne \lamTwo | \DeKhat | \left[ s_K \brOne^T \frac{\GDeK}{| \DeK |} \brOne + s_K^{-1} \brTwo^T \frac{\GDeK}{| \DeK |} \brTwo \right] ,
\end{equation}
where $| \DeK | {=} \lamOne \lamTwo | \DeKhat |$, $\DeKhat {=} T_K^{-1}(\DeK)$ being the pullback of the patch $\DeK$ via the transformation $T_K$. More precisely, the information on the area of the element is encapsulated in the term $\lamOne \lamTwo | \DeKhat |$, whereas the remaining terms strictly depend on the anisotropic information of the element via $\{\brOne, \brTwo \}$ and the aspect ratio $s_K$.

\subsection{Construction of the anisotropically adapted mesh}
\label{sc:ConstructMesh}

The error estimate in equation~\eqref{eq:estimateLoc} provides information on the regions where additional accuracy is required to describe the segmentation contour. In order to improve such a description, a metric-based procedure is employed to construct an adapted mesh~\cite{George-BGHLS-96,Frey-FG-07}.

Starting from the current mesh $\Th$, the local information provided by the error estimate $\eta_K$ is exploited in a predictive way to construct a symmetric positive definite matrix $\Metric {\in} \RR^{2 \times 2}$, known as \emph{metric}, containing all the geometric information of the target mesh. From a practical viewpoint, the metric is approximated on the mesh $\Th$ by means of a piecewise constant function such that
\begin{equation}
\label{eq:metricLoc}
\Metric |_K = \RK^T \LK^{-2} \RK \quad \forall  K \in \Th.
\end{equation}

It can be verified, see~\cite{Perotto-MP-06}, that a new metric providing a mesh with minimum cardinality for a user-defined  accuracy on the error is obtained via an iterative procedure minimising, for each element, the functional
\begin{equation}
\label{eq:meshMinPb}
\Jfunc(s_K, \brOne, \brTwo) =  s_K \brOne^T \frac{\GDeK}{| \DeK |} \brOne + s_K^{-1} \brTwo^T \frac{\GDeK}{| \DeK |} \brTwo ,
\end{equation}
under the constraints $s_K {\geq} 1$ and $\bm{r}_{i,K} {\cdot} \bm{r}_{j,K} {=} \delta_{ij}$, $\delta_{ij}$ being the Kronecker delta.
The solution of this optimisation problem is provided by the following result (see~\cite{Perotto-MP-06} for the proof):

\begin{proposition}
Let $\{ \theta_i, \bt_i \}, \ i {=} 1,2$ be the eigenvalues and eigenvectors of the matrix $\GDeK / | \DeK |$, with $\bt_1$ and $\bt_2$ orthonormal vectors and $\theta_1 \geq \theta_2 > 0$. The values $(s_K^\star, \brOne^\star, \brTwo^\star)$ minimising the functional in equation~\eqref{eq:meshMinPb} are given by
\begin{equation}
\label{eq:optimalValues}
s_K^\star  = (\theta_1 / \theta_2)^{1/2} , \quad \brOne^\star = \bt_2 , \quad \brTwo^\star = \bt_1 .
\end{equation}
\end{proposition}

To fully describe the metric $\Metric^\star$ associated with the adapted mesh, the values of the optimal eigenvalues $\lamOne^\star$ and $\lamTwo^\star$, now merged in $s_K^\star$, need to be individually identified for each element $K$. It follows that
%
\begin{subequations}
\label{eq:lambdaOpt}
\begin{align}
\lamOne^\star &= \theta_2^{-1/2} \left(\frac{\tau^\star \left\| \grad \phi_h \right\|_{\DeK}^2 }{2 | \Khat |} \right)^{1/2} , \\
\lamTwo^\star &= \theta_1^{-1/2} \left(\frac{\tau^\star \left\| \grad \phi_h \right\|_{\DeK}^2 }{2 | \Khat |} \right)^{1/2} ,
\end{align}
\end{subequations}
where $\Khat {=} T_K^{-1}(K)$ is the pullback of element $K$ via the map $T_K$ and $\tau^\star$ is a user-defined accuracy. It is worth noticing that the optimal eigenvalues in equation~\eqref{eq:lambdaOpt} are obtained imposing a local tolerance on the error in element $K$. More precisely, the local tolerance $\tau^\star \left\| \grad \phi_h \right\|_{\DeK}^2 $ prescribes a target precision in element $K$, as a relative measure of the average Euclidean norm of the gradient of the level-set solution on the patch $\DeK$, weighted by the area of the patch elements, that is,
\begin{equation}
\label{eq:refMeasure}
\left\| \grad \phi_h \right\|_{\DeK}^2 := \frac{1}{| \DeK |} \sum_{T \in \DeK} | T | \left\| \grad \phi_h|_T \right\|^2 .
\end{equation}
Hence, the optimal piecewise constant metric $\Metric^\star$ is provided by the pairs $\{ \lambda_{i,K}^\star, \bm{r}_{i,K}^\star \}, \ i {=} 1,2$ for all mesh elements $K {\in} \Th$. The corresponding adapted mesh is obtained by means of the \texttt{adaptmesh} routine available in \texttt{FreeFem++}~\cite{Hecht-12}.

\begin{remark}[Maximum stretching]
According to the above procedure, the optimal eigenvalues $\lamOne^\star$ and $\lamTwo^\star$ are computed starting from the value of the optimal stretching factor $s_K^\star$. Since the level-set function $\phi$ is discretised by means of a piecewise constant approximation, an excessive stretching of the mesh elements could be responsible for a loss of accuracy of the approximation in the direction parallel to the interface. Hence, from a practical viewpoint, the mesh adaptation procedure introduces a user-defined upper bound for $s_K^\star$. In the  simulations in section~\ref{sc:Simulations}, the maximum admissible value for the stretching factor is fixed to $1,000$.
\end{remark}

\begin{remark}[Metric relaxation]
In order to avoid abrupt changes in the element size during the adaptation procedure, a relaxation step is introduced when defining the new metric $\MetricOmega^\star$. Hence, this is obtained as a linear combination of the computed target metric $\Metric^\star$ with the previous metric $\MetricOld$, namely
\begin{equation}
\label{eq:metricComb}
\MetricOmega^\star = \omega \Metric^\star + (1-\omega) \MetricOld , \quad \omega \in [0,1] .
\end{equation}
\end{remark}

\subsection{Split-adapt Bregman algorithm for image segmentation}
\label{sc:SplitAdaptBreg}

The split Bregman method presented in algorithm~\ref{alg:splitBregman} is now enriched with the anisotropic mesh adaptation procedure just detailed. Indeed, the discrete level-set function $\phi_h$ exhibits strong gradients across the contour of the region to be segmented. An anisotropically adapted mesh is thus ideal to capture such directional features, providing an accurate  description of the interface between the background and the foreground of the image, with a reduced number of degrees of freedom. The resulting strategy is described in algorithm~\ref{alg:splitAdaptBregman}.

\begin{algorithm}
\caption{Split-adapt Bregman algorithm for image segmentation}\label{alg:splitAdaptBregman}
\begin{algorithmic}[1]
\REQUIRE{Tolerance $\eta^\star$ for the stopping criterion and initial contour $\phi_0$. 
For adaptation: target accuracy $\tau^\star$, number $\nRec$ of Bregman iterations between two adaptation steps and metric relaxation parameter $\omega$.}
\STATE{Construct the initial computational mesh $\Th^0$.}
\STATE{Set $k {=} 0$, $\phi_h^0 {=} \phi_0$, $\bd_h^0 {=} \grad \phi_0$ and $\bb_h^0 {=} \bm{0}$.}
\WHILE{$\left\| \phi_h^{k+1} - \phi_h^k \right\|_2 > \eta^\star \left\| \phi_h^k \right\|_2$}
\STATE{Compute the contour $\phi_h^{k+1}$ solving equation~\eqref{eq:evolutionPhiFE}.}
\IF{$k \operatorname{mod} \nRec = 0$}
\STATE{Construct the gradient reconstruction $\PZZ \phi_h^{k+1}$ using~\eqref{eq:recoveryConst}.}
\STATE{Use $\tau^\star$ to estimate the eigenpairs $\{ \lambda_{i,K}^\star, \bm{r}_{i,K}^\star \} , i {=} 1,2$ in~\eqref{eq:optimalValues}-\eqref{eq:lambdaOpt}.}
\STATE{Construct the relaxed metric $\MetricOmega^\star$ via~\eqref{eq:metricComb}.}
\STATE{Generate the mesh $\Th^{k+1}$ using $\MetricOmega^\star$ and the \texttt{adaptmesh} routine.}
\STATE{Interpolate $\phi_h^{k+1}$, $\bd_h^k$ and $\bb_h^k$ on the new mesh $\Th^{k+1}$.}
\STATE{Reduce the value of the target accuracy $\tau^\star \gets \tau^\star / 2$.}
\ENDIF
\STATE{Compute the auxiliary variable $\bd_h^{k+1}$ via the shrinkage operator~\eqref{eq:shrinkD}.}
\STATE{Compute the Bregman update $\bb_h^{k+1}$ using equation~\eqref{eq:splitBregman-C}.}
\STATE{$k \gets k+1$.}
\ENDWHILE
\ENSURE{Boundary $\phi_h^k$ of the segmented region and final adapted mesh $\Th^k$.}
\end{algorithmic}
\end{algorithm}

Following the structure of split Bregman, the split-adapt Bregman algorithm first solves step A to compute an approximation $\phi_h^{k+1}$ of the level-set function after fixing the last available values for $\bd_h$ and $\bb_h$ at iteration $k$ (Algorithm~\ref{alg:splitAdaptBregman} - Step 4). Then, every $\nRec$ Bregman iterations (Algorithm~\ref{alg:splitAdaptBregman} - Step 5), the mesh adaptation procedure is performed: 
\begin{enumerate}[label=(\roman*)]
\item the reconstructed gradient $\PZZ \phi_h$ (Algorithm~\ref{alg:splitAdaptBregman} - Step 6) is employed to estimate the eigenpairs of the target metric (Algorithm~\ref{alg:splitAdaptBregman} - Step 7);
\item the relaxed metric $\MetricOmega^\star$ is constructed (Algorithm~\ref{alg:splitAdaptBregman} - Step 8) and employed to construct the associated adapted mesh $\Th^{k+1}$, avoiding any abrupt change in the element size (Algorithm~\ref{alg:splitAdaptBregman} - Step 9);
\item the last computed finite element approximations are interpolated on the new mesh $\Th^{k+1}$ (Algorithm~\ref{alg:splitAdaptBregman} - Step 10);
\item the target accuracy $\tau^\star$ of the adaptation routine is halved (Algorithm~\ref{alg:splitAdaptBregman} - Step 11); notice that this step is performed for a maximum of 5 times, then the value of $\tau^\star$ is frozen.
\end{enumerate}
Finally, steps B and C perform the shrinkage operation and the Bregman update to determine $\bd_h^{k+1}$ and $\bb_h^{k+1}$, respectively (Algorithm~\ref{alg:splitAdaptBregman} - Steps 13-14), before the iteration counter $k$ is updated (Algorithm~\ref{alg:splitAdaptBregman} - Step 15) and a new iteration is carried out. As for the split Bregman method, the algorithm stops when the relative variation of the level-set function between  two consecutive iterations is below a user-defined tolerance $\eta^\star$ (Algorithm~\ref{alg:splitAdaptBregman} - Step 3). 

Similarly to the split Bregman algorithm, also the split-adapt Bregman strategy is suitable for any image segmentation model that can be formulated as in equation~\eqref{eq:minL1phi}.

\section{Numerical experiments}
\label{sc:Simulations}

In this section, the performance of the split-adapt Bregman algorithm is analysed by means of a set of numerical experiments using both synthetic and real images. The proposed approach outperforms the standard split Bregman strategy by enhancing the quality of the segmented region at a limited extra cost. Moreover, the split-adapt Bregman algorithm is shown to be robust to different  sources of noise, including Gaussian, salt and pepper and speckle noise. Finally, the method is applied to the segmentation of a challenging medical image.

\subsection{Validation of the split-adapt Bregman algorithm}
\label{sc:Validation}

\begin{figure}[!htb]
	\centering
	\subfigure[Image \label{fig:validationImg}]{\includegraphics[width=0.32\textwidth]{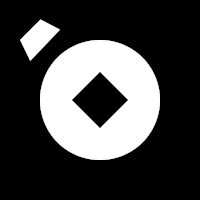}}
	\subfigure[RSFE - Split Bregman \label{fig:validationRSFE}]{\includegraphics[width=0.32\textwidth]{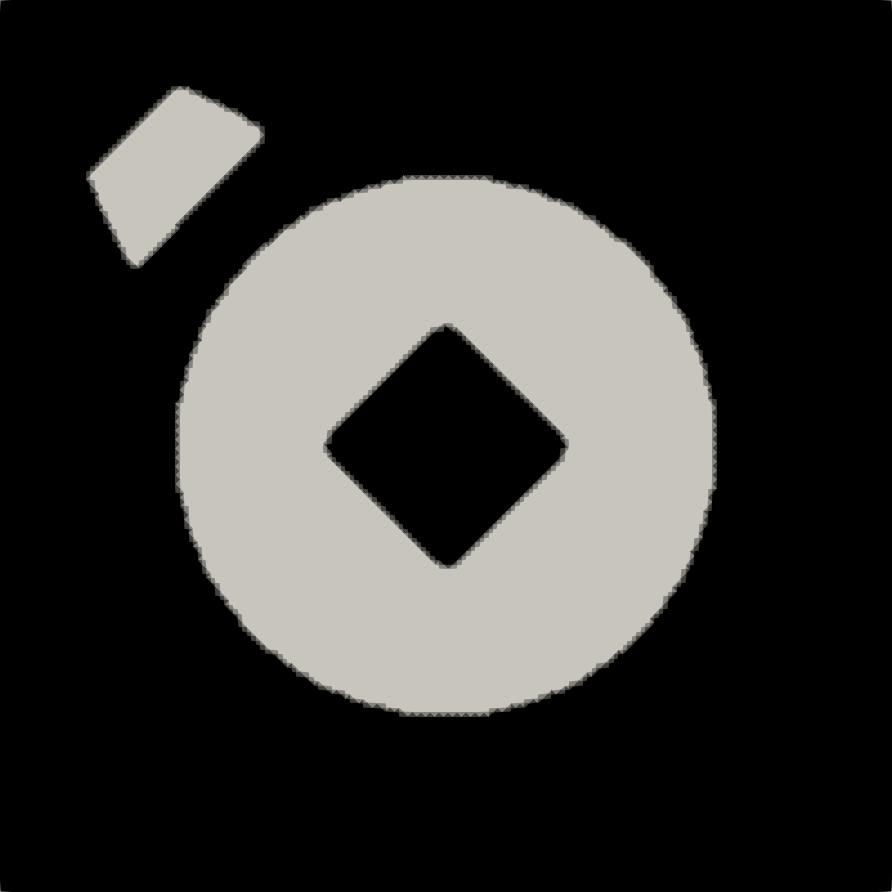}}	
	\subfigure[Bayesian - Split Bregman \label{fig:validationBayes}]{\includegraphics[width=0.32\textwidth]{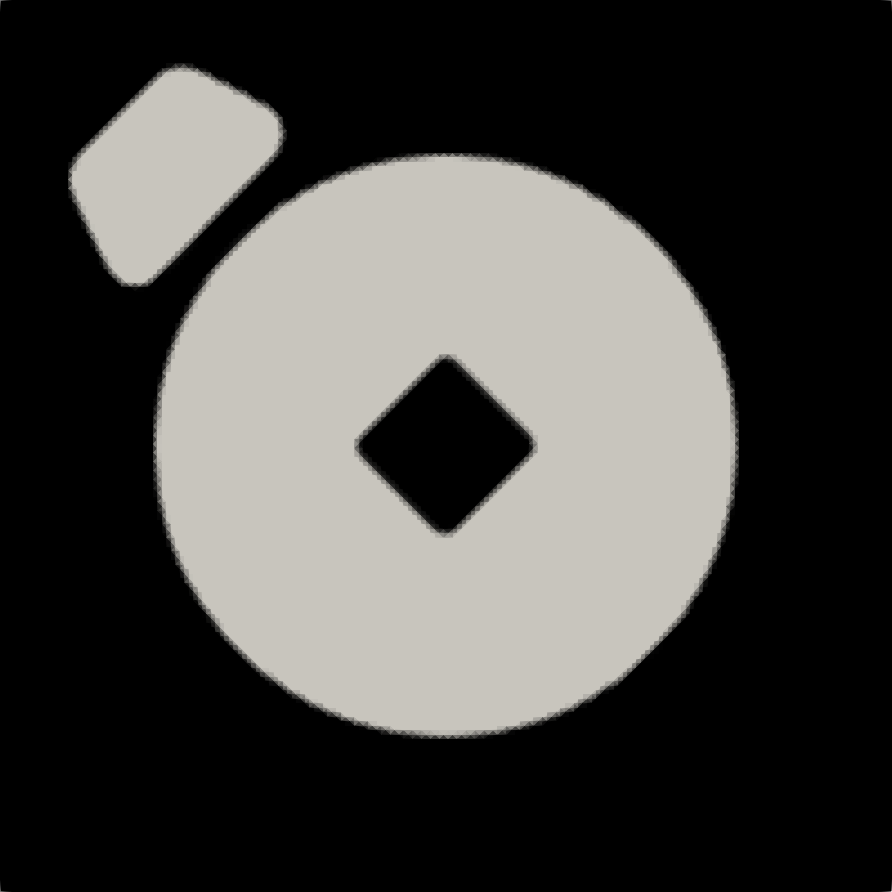}}
		
	\subfigure[Initial guess \label{fig:validationInit}]{\includegraphics[width=0.32\textwidth]{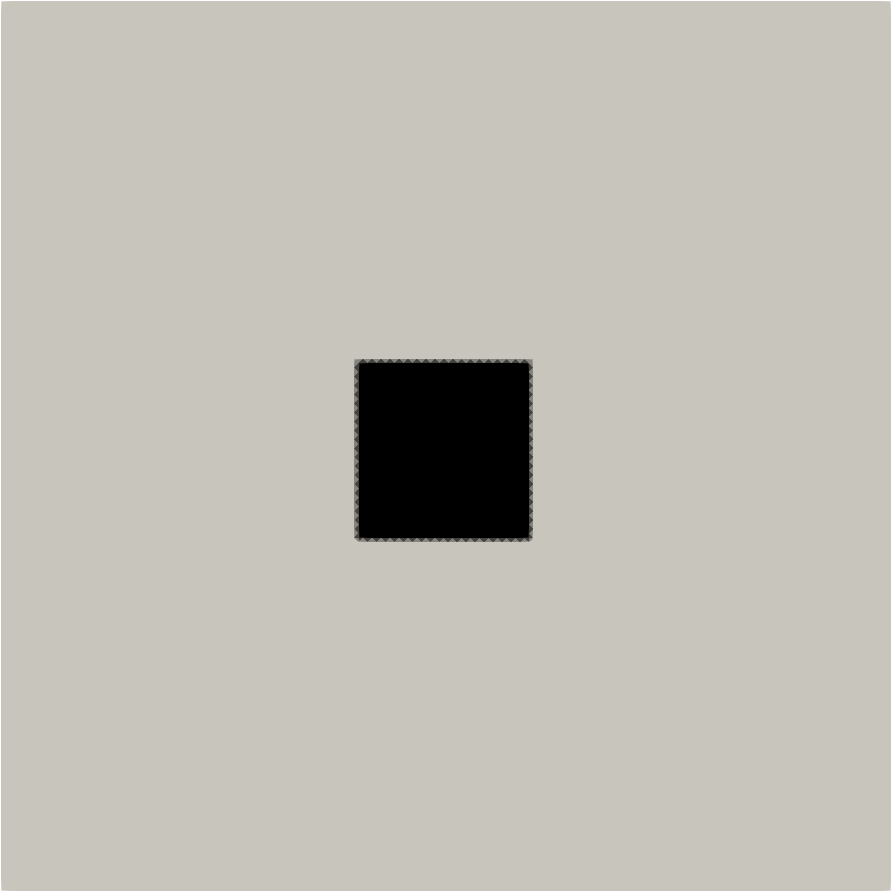}}
	\subfigure[Detail of~\ref{fig:validationRSFE} \label{fig:validationRSFE_zoom}]{\includegraphics[width=0.32\textwidth]{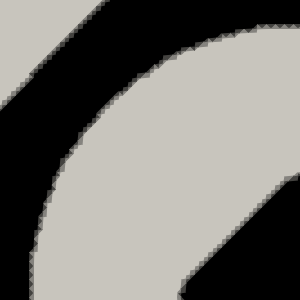}}	
	\subfigure[Detail of~\ref{fig:validationBayes} \label{fig:validationBayes_zoom}]{\includegraphics[width=0.32\textwidth]{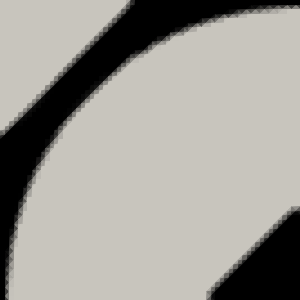}}	
	\caption{Synthetic image with homogeneous regions - (a) Image and (d) initial guess for the segmentation. Final segmented results on a uniform structured mesh using the standard split Bregman algorithm with (b) the RSFE model and (c) the Bayesian approach. (e-f) Details of the segmentation outcomes.}
\label{fig:validation}
\end{figure}

A validation of the split-adapt Bregman method is performed using a $200 \times 200$ pixels synthetic image with homogeneous regions, see figure~\ref{fig:validationImg}. In this context, the region-based Bayesian model presented in section~\ref{sc:ImgSeg} is compared with the classical RSFE model described in~\ref{sc:appRSFE}.

The problem is first solved using the standard split Bregman algorithm on a domain discretised with a uniform structured mesh of 79,202 triangular elements, with $\hMin {=} 1$ and $\hMax {=} 1.4$. The initial guess for the segmentation is displayed in figure~\ref{fig:validationInit}. Both the RSFE model and the Bayesian strategy achieve the required tolerance of $\eta^\star {=} 0.5 \times 10^{-2}$. More precisely, the RSFE model requires 12 iterations, whereas the Bayesian strategy converges in 8 iterations. The two approaches are able to subdivide the regions within the image, also identifying non-connected areas as it is evident in figures~\ref{fig:validationRSFE} and~\ref{fig:validationBayes}. In this example, the RSFE model provides more accuracy than the Bayesian approach, as shown by the perfect match of the size of the segmented regions in figure~\ref{fig:validationRSFE} with respect to the original image in figure~\ref{fig:validationImg}, whereas the result of the Bayesian approach in figure~\ref{fig:validationBayes} appears slightly overdiffusive, partially smoothing the boundaries and overgrowing the segmented regions.

The description of the region boundaries is clearly affected by the quality of the underlying computational mesh. Figures~\ref{fig:validationRSFE_zoom} and~\ref{fig:validationBayes_zoom} display a detail of the segmented contours provided by the split Bregman method applied to the RSFE and the Bayesian model, respectively. Although the latter strategy provides a more regular approximation of the contour reducing the jagged effect which characterises the RSFE solution, the smoothness of the interface is  still severely affected by the lack of resolution of the computational grid. 

\begin{figure}[!htb]
	\centering
	\subfigure[RSFE - Split-adapt Bregman \label{fig:validationRSFE_ad}]{\includegraphics[width=0.45\textwidth]{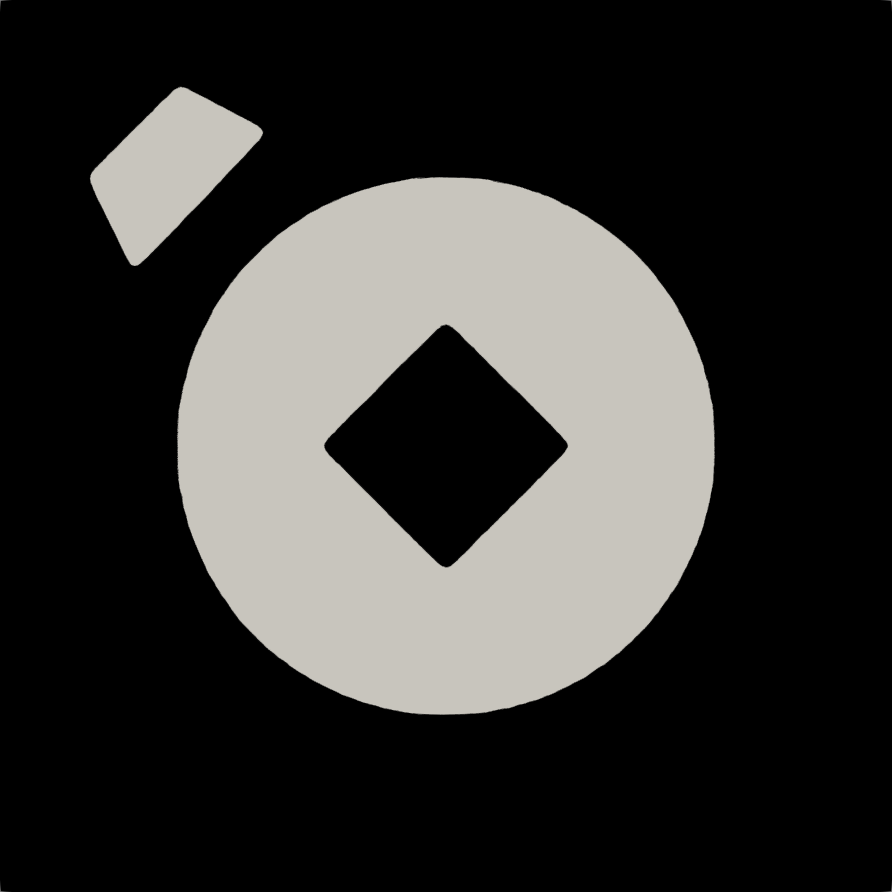}}	
	\subfigure[Bayesian - Split-adapt Bregman \label{fig:validationBayes_ad}]{\includegraphics[width=0.45\textwidth]{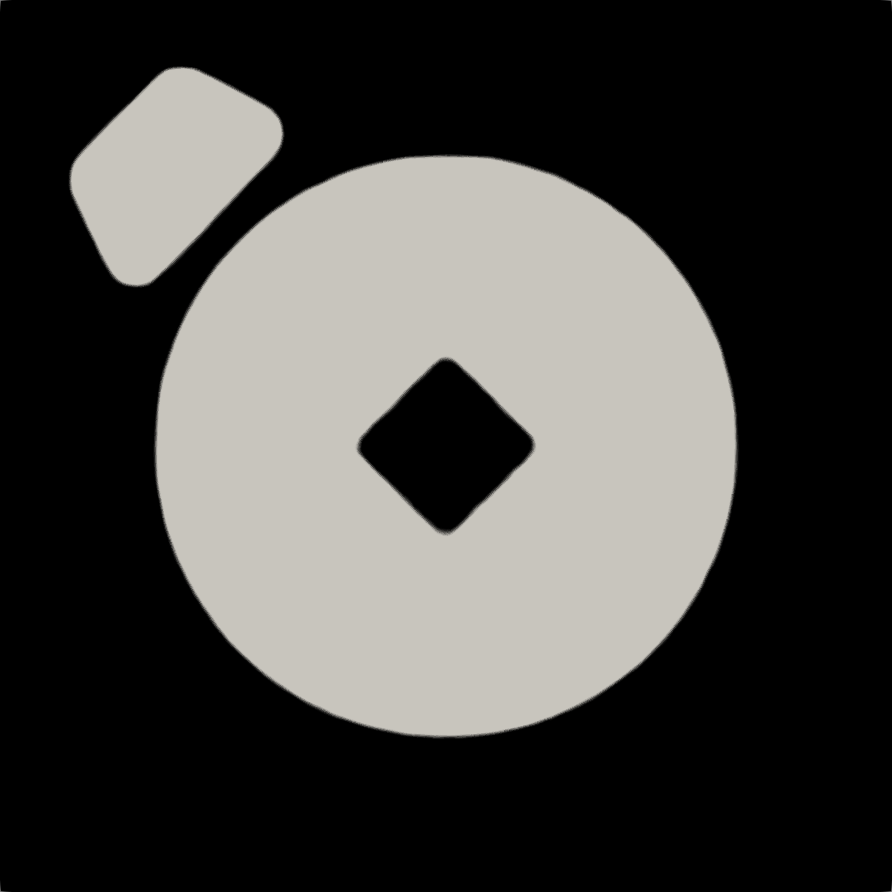}}
		
	\subfigure[Detail of~\ref{fig:validationRSFE_ad}]{\includegraphics[width=0.45\textwidth]{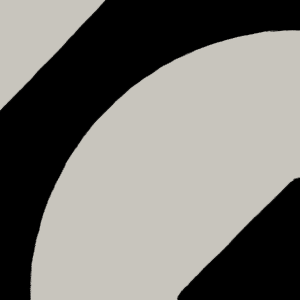}}	
	\subfigure[Detail of~\ref{fig:validationBayes_ad}]{\includegraphics[width=0.45\textwidth]{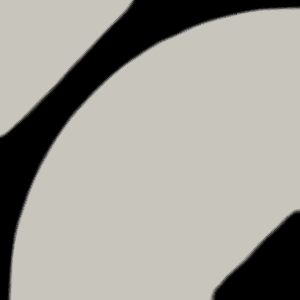}}	
	\caption{Synthetic image with homogeneous regions - Final segmented results on an adapted mesh using the split-adapt Bregman algorithm with (a) the RSFE model and (b) the Bayesian approach. (c-d) Details of the segmentation outcomes.}
\label{fig:validationAdapt}
\end{figure}

The numerical test is thus repeated using the split-adapt Bregman algorithm with the same value of the tolerance $\eta^\star$ and the same initial guess. This strategy performs one anisotropic mesh adaptation step every $\nRec {=} 3$ iterations of the Bregman solver. The tolerance for the error estimate is initially set to $\tau^\star {=} 0.5$ and it is successively halved after each run of the adaptation routine. In addition, the relaxation parameter $\omega {=} 0.9$ is selected to construct  the new metric. The split-adapt Bregman algorithm converges after 12 iterations with 4 mesh adaptation steps for the RSFE model, whereas 7 iterations with 2 mesh adaptation steps are required for the Bayesian model. The outcome of the segmentation, provided in figure~\ref{fig:validationAdapt}, highlights the superiority of the split-adapt Bregman algorithm with respect to the standard split Bregman approach, in terms of accuracy of the interface description. As previously observed, for this example the segmentation yielded by the RSFE model provides excellent accuracy in both the shape and the size of the segmented regions, whereas the Bayesian approach slightly overgrows the segmented regions.

\begin{figure}[!htb]
	\centering
	\subfigure[RSFE - Final mesh \label{fig:validationRSFE_mesh}]{\includegraphics[width=0.45\textwidth]{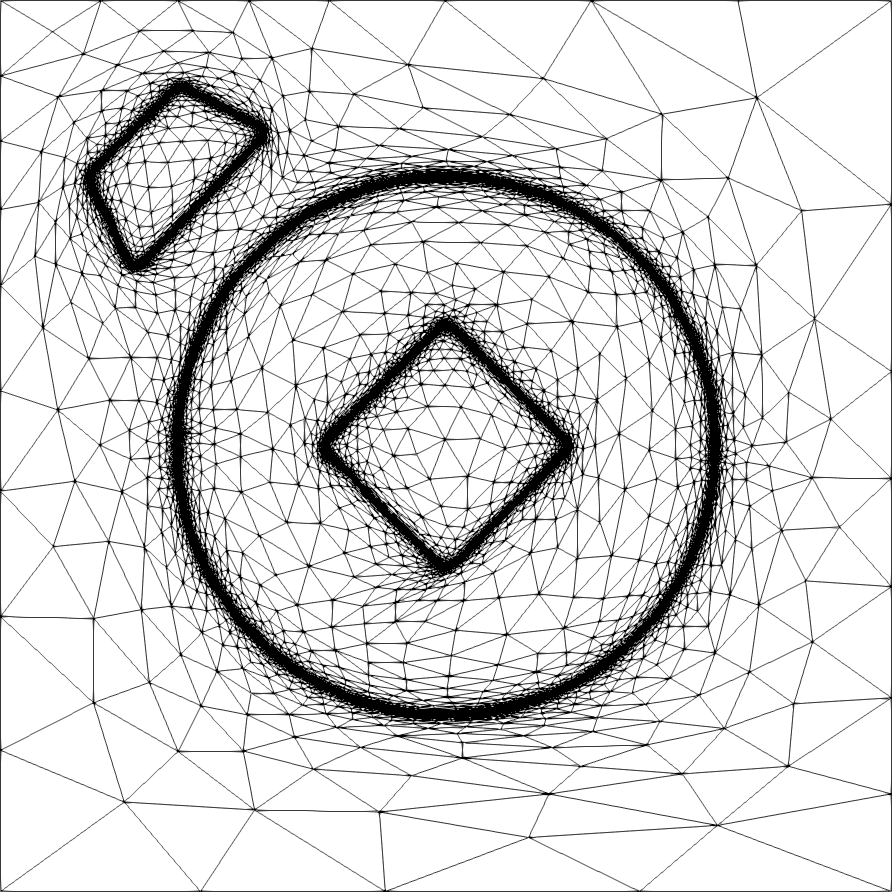}}	
	\subfigure[Bayesian - Final mesh \label{fig:validationBayes_mesh}]{\includegraphics[width=0.45\textwidth]{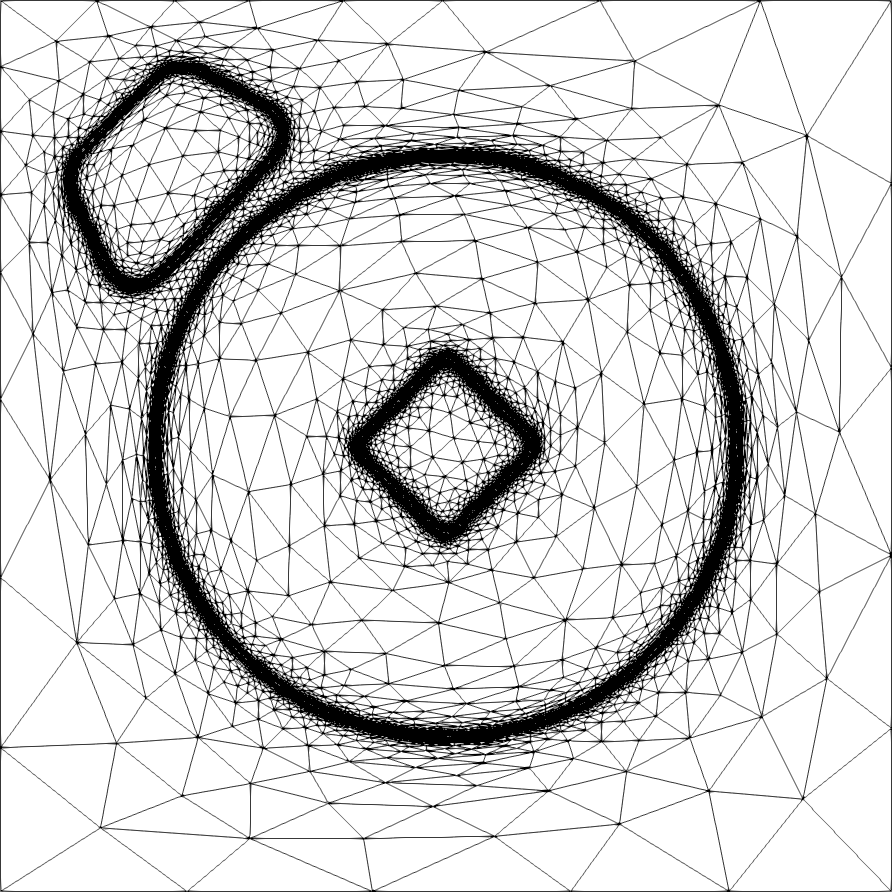}}
		
	\subfigure[Detail of~\ref{fig:validationRSFE_mesh} \label{fig:validationRSFE_mesh_zoom}]{\includegraphics[width=0.45\textwidth]{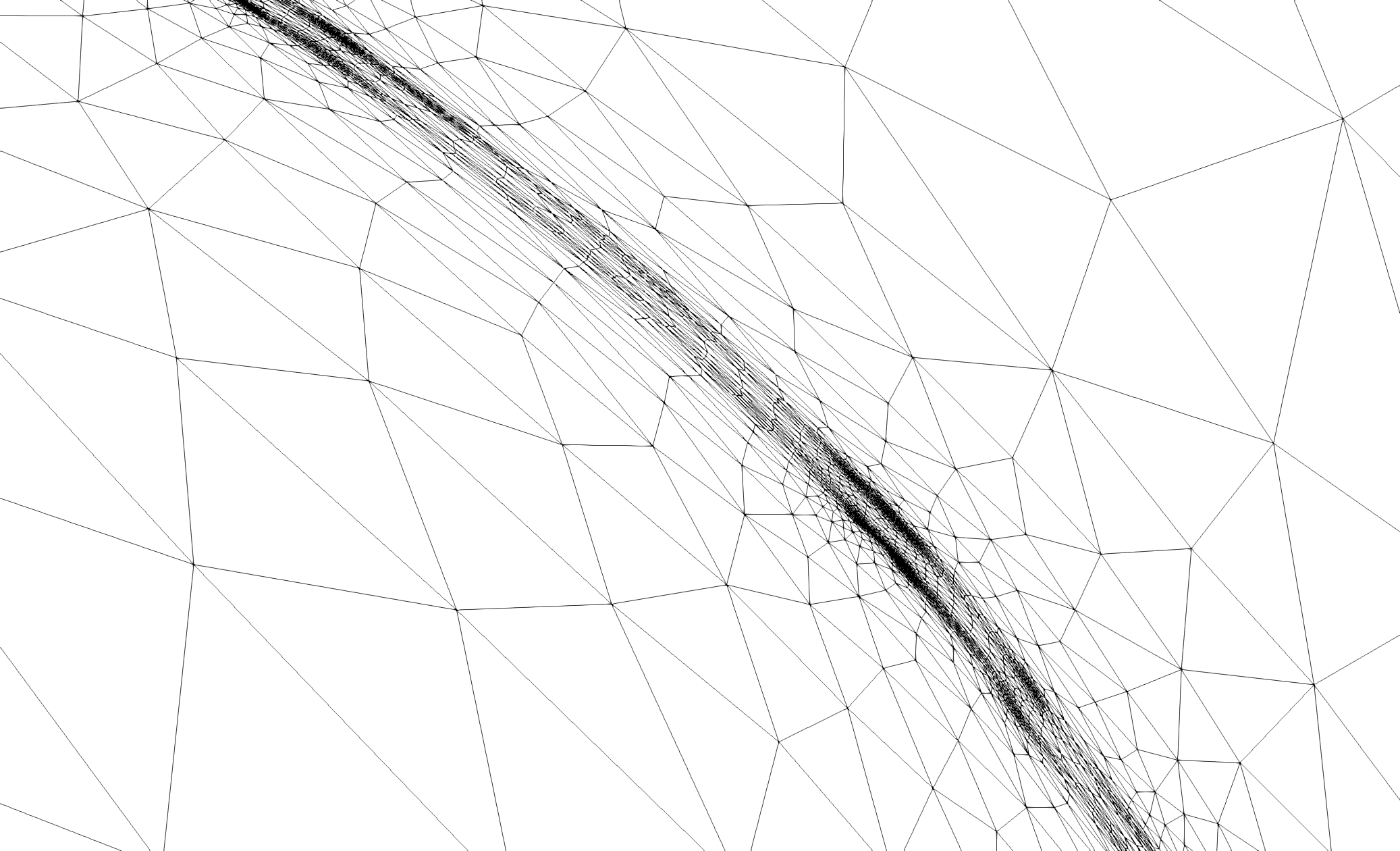}}	
	\subfigure[Detail of~\ref{fig:validationBayes_mesh} \label{fig:validationBayes_mesh_zoom}]{\includegraphics[width=0.45\textwidth]{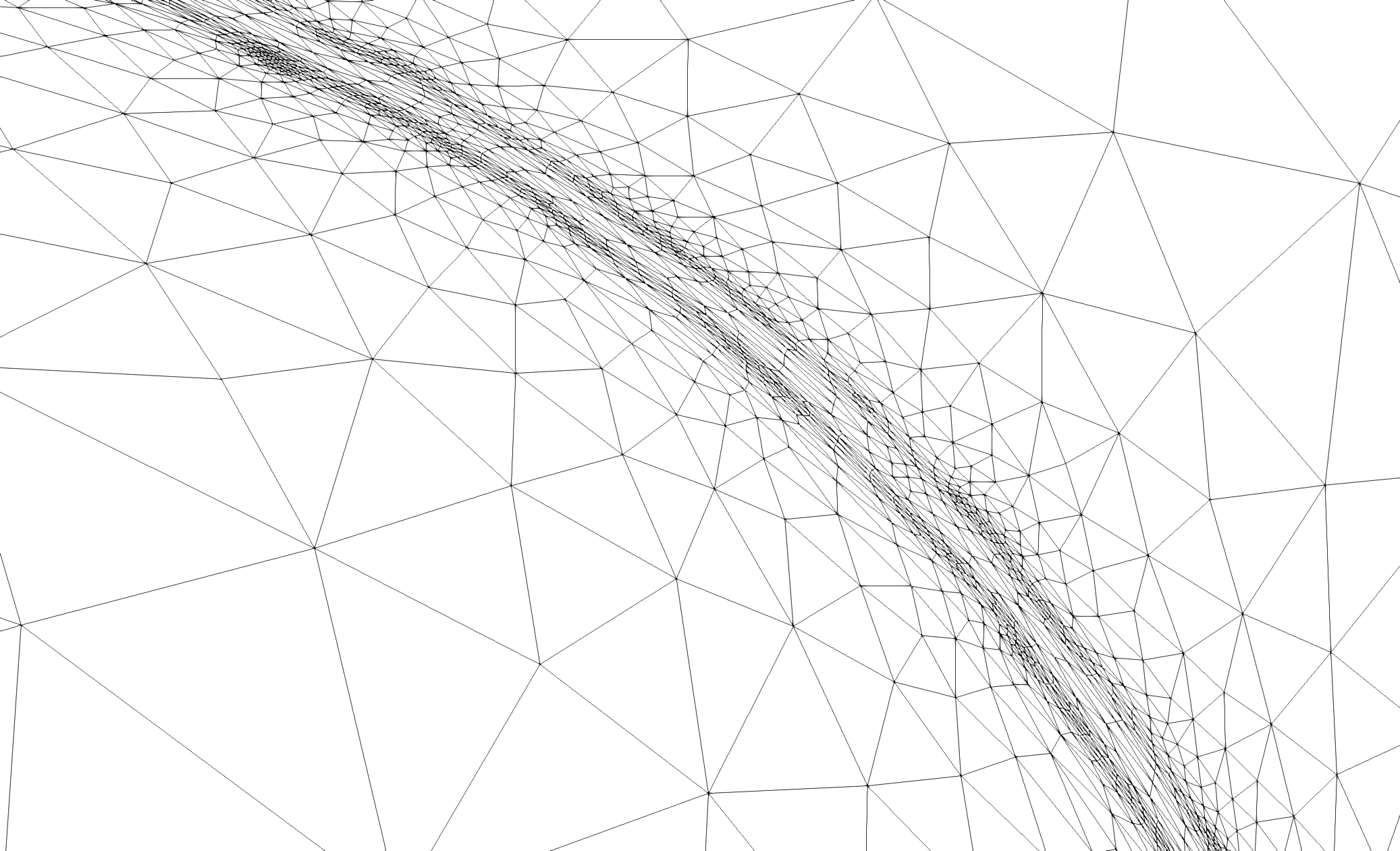}}	
	\caption{Synthetic image with homogeneous regions - Final adapted mesh using the split-adapt Bregman algorithm with (a) the RSFE model and (b) the Bayesian approach. (c-d) Details of the final meshes.}
\label{fig:validationMesh}
\end{figure}

Although the split-adapt Bregman algorithm entails an extra cost associated with the computation of the new metric and the adaptation step, the resulting meshes feature a reduced number of elements with respect to the uniform structured ones. Hence, a significant computational saving characterises the Bregman steps. More precisely, the final mesh for the RSFE model, see figure~\ref{fig:validationRSFE_mesh}, contains 29,281 triangular elements, with $\hMin {=} 0.27 \times 10^{-1}$, $\hMax {=} 35.69$ and a maximum elemental stretching factor $s_K {=} 844.99$. Similarly, figure~\ref{fig:validationBayes_mesh} displays the final mesh for the Bayesian approach featuring 55,614 triangles, with $\hMin {=} 0.1$, $\hMax {=} 18.32$ and maximum aspect ratio $s_K {=} 822.48$. 
It is worth noticing that, although the maximum stretching factor of the elements is comparable in the two cases, the adapted mesh in figure~\ref{fig:validationBayes_mesh} maintains regions with isotropic elements, compare, e.g., the detail in figure~\ref{fig:validationBayes_mesh_zoom} with the corresponding one, obtained using the RSFE functional, in figure~\ref{fig:validationRSFE_mesh_zoom}. Such a result can be ascribed to the reduced number of adaptation steps performed for the Bayesian model and to the relaxation of the initial computational grid, resulting in a final mesh with a larger number of elements with respect to the RSFE model. Nonetheless, the overall reduced number of iterations of the optimisation algorithm (i.e., 7 versus 12) still makes the split-adapt Bregman method for the Bayesian functional competitive from the computational viewpoint.

\begin{remark}[Application to other segmentation models]
The above simulations show that the split-adapt Bregman algorithm works seamlessly both with the RSFE and the Bayesian model, without requiring any interaction with the user. It is worth mentioning that the applicability of the split-adapt Bregman methodology can be generalised to any region-based segmentation model that can be expressed according to equation~\eqref{eq:minL1phi} and that can be solved using the standard split Bregman approach.
\end{remark}

\begin{remark}[Frequency of mesh adaptation]
A key aspect to control the overall computational cost of the split-adapt Bregman algorithm is the selection of the number $\nRec$ of iterations of the optimisation method between two adaptation steps. Numerical experiments, not reported here for brevity, were performed to test different values of this parameter. More precisely, for $\nRec {=} 1$, adaptation is performed after every Bregman iteration, introducing a significant additional cost into the algorithm for the generation of the new mesh and the interpolation of the previously computed finite element approximations. This approach was disregarded for its excessive computational burden. For large values of $\nRec$, the adapted mesh is frozen for several iterations of the optimisation algorithm. In this case, after the initial iterations, the evolving interface is not properly tracked since the adapted mesh is excessively refined in some regions of the domain, while the discretisation of other areas is too coarse, resulting in a reduced accuracy of the representation of the boundary.
Hence, the value $\nRec {=} 3$ was selected for all the simulations in the present work, as a trade-off between computational cost and accuracy of the split-adapt Bregman algorithm.
\end{remark}

\subsection{Treatment of image spatial information}
\label{sc:Noise}

\begin{figure}[!htb]
	\centering
	\subfigure[Image \label{fig:noiseImg}]{\includegraphics[width=0.32\textwidth]{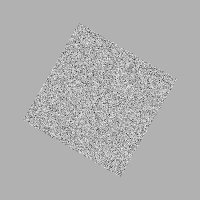}}
	\subfigure[RSFE - Split Bregman \label{fig:noiseRSFE}]{\includegraphics[width=0.32\textwidth]{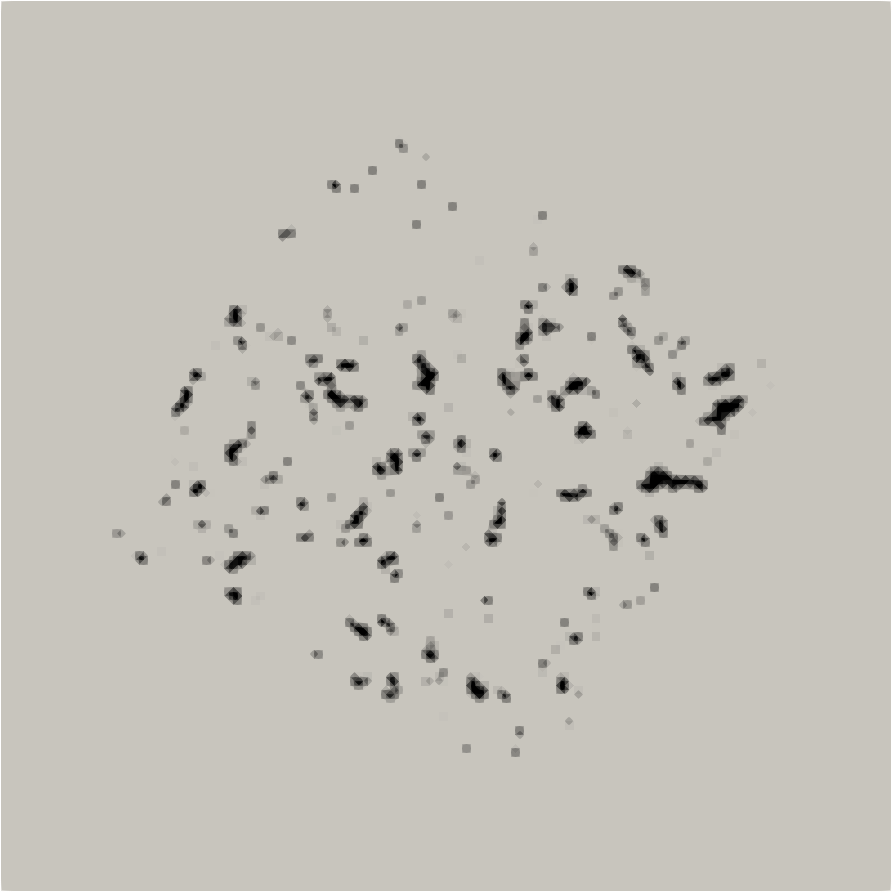}}	
	\subfigure[Bayesian - Split Bregman \label{fig:noiseBayes}]{\includegraphics[width=0.32\textwidth]{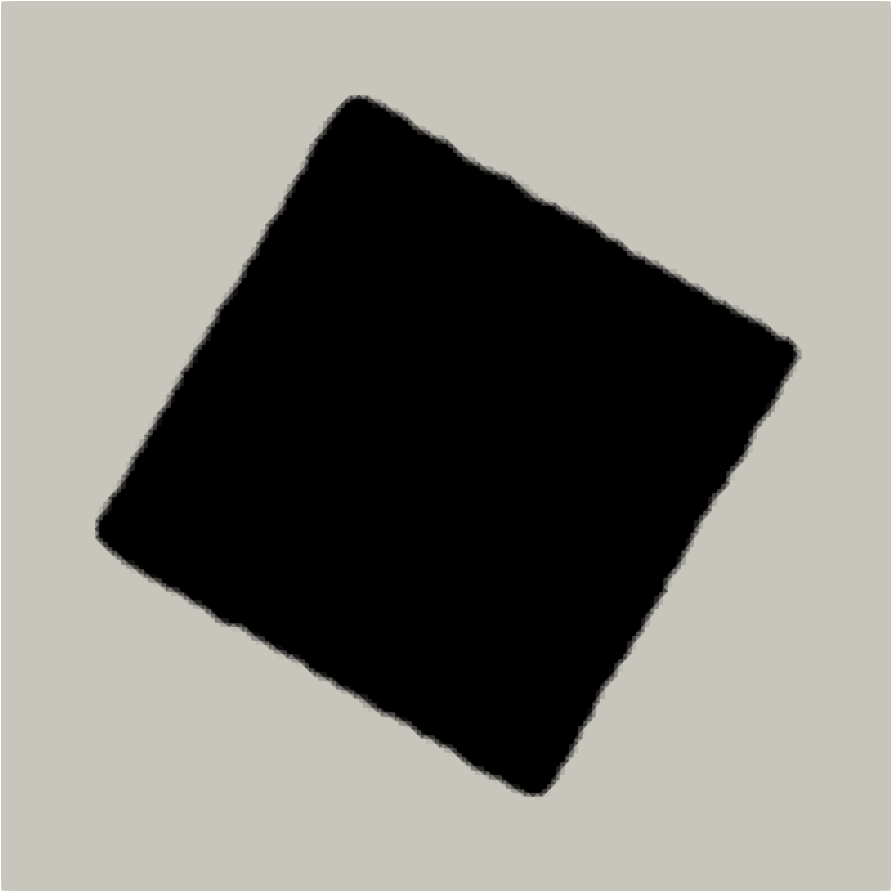}}
		
	\caption{Synthetic image with spatial information - (a) Image. Final segmented results on a uniform structured mesh using the split Bregman algorithm with (b) the RSFE model after 500 iterations (no  convergence is attained) and (c) the Bayesian approach after 12 iterations.}
\label{fig:noise}
\end{figure}

For the second test, a $200  \times 200$ pixels synthetic image with mean value of the RGB intensity equal to $[177,177,177]$ and non-uniform variance is considered, see figure~\ref{fig:noiseImg}. The variance of the background is set to $0.05$, whereas a square region in the centre of the domain is characterised by variance $0.4$, thus introducing a spatial information of the pixel inhomogeneity  which represents a challenge for image segmentation algorithms.

\begin{figure}[!htb]
	\centering
	\subfigure[External region]{\includegraphics[width=0.49\textwidth]{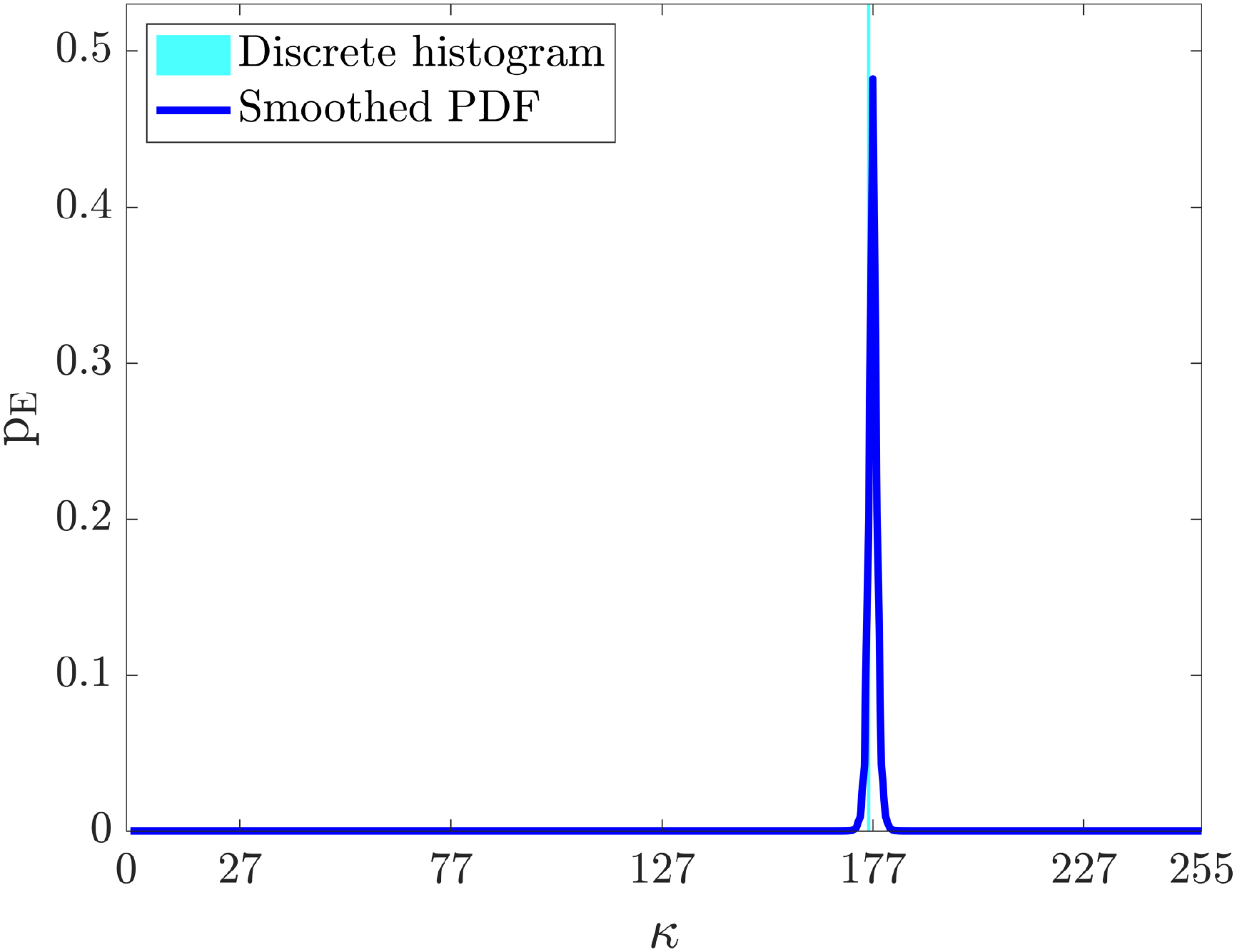}}
	\subfigure[Internal region]{\includegraphics[width=0.49\textwidth]{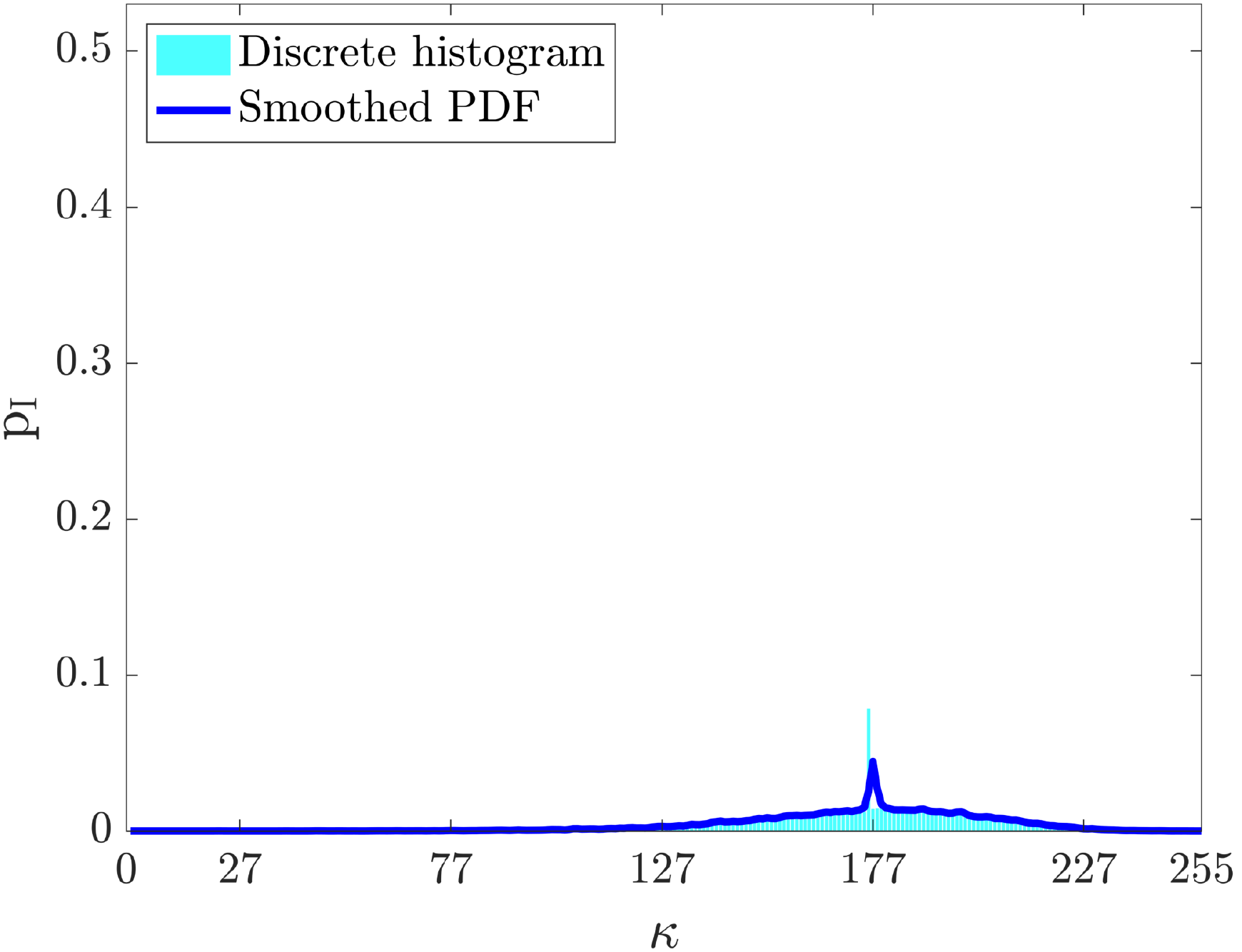}}	
		
	\caption{Synthetic image with spatial information - Estimated probability density functions of the pixel intensity $\kappa$ for the (a) external and (b) internal region.}
\label{fig:noiseDensity}
\end{figure}

As for the case in the previous section, the segmentation problem is first solved using the standard split Bregman method to minimise the RSFE and the Bayesian functional. The domain is discretised using a uniform structured grid of 79,202 triangles, with $\hMin {=} 1$ and $\hMax {=} 1.4$. The initial guess for the segmentation is the same as in figure~\ref{fig:validationInit}, while the tolerance for the stopping criterion is set to $\eta^\star {=} 0.5 \times 10^{-3}$. The split Bregman algorithm applied to the RSFE model is unable to converge. The solution obtained after 500 iterations is shown in figure~\ref{fig:noiseRSFE}, highlighting the failure of the RSFE functional in correctly identifying the two regions in the image. More precisely, the RSFE model only considers a small region surrounding each pixel through the convolution with a Gaussian kernel, see equation~\eqref{eq:fIntExt}, and it does not fully exploit the image spatial information. On the contrary, the Bayesian model embeds such a spatial information into the PDFs of the pixel intensity in each region, by performing the nonparametric kernel density estimate in equation~\eqref{eq:pIntExt}. The estimated PDFs obtained from the smoothing of the discrete histograms of the pixel intensity $\kappa$ are displayed in figure~\ref{fig:noiseDensity}. Both probability distributions are centred on the value $\kappa {=} 177$. Moreover, the external region features a limited variation of the intensity, whereas the internal one has a larger variance, coherently with the original image in figure~\ref{fig:noiseImg}. The solution provided by the split Bregman algorithm with the Bayesian functional converges in 12 iterations and the resulting segmentation accurately identifies the two regions in the image as it can be appreciated in figure~\ref{fig:noiseBayes}. 

Although the segmentation result in figure~\ref{fig:noiseBayes} is able to discriminate the presence of two regions by exploiting the  spatial information of the image, the quality of the description of the interface is still affected by the poor resolution of the underlying computational grid. The split-adapt Bregman algorithm is thus applied to the Bayesian model to improve the quality of the segmentation. \hl{The complete setup of the input parameters required by the method is detailed in~\ref{sc:appParam}.}

\begin{figure}[!htb]
	\centering
	\subfigure[Final segmentation  \label{fig:noiseBayes_ad}]{\includegraphics[width=0.45\textwidth]{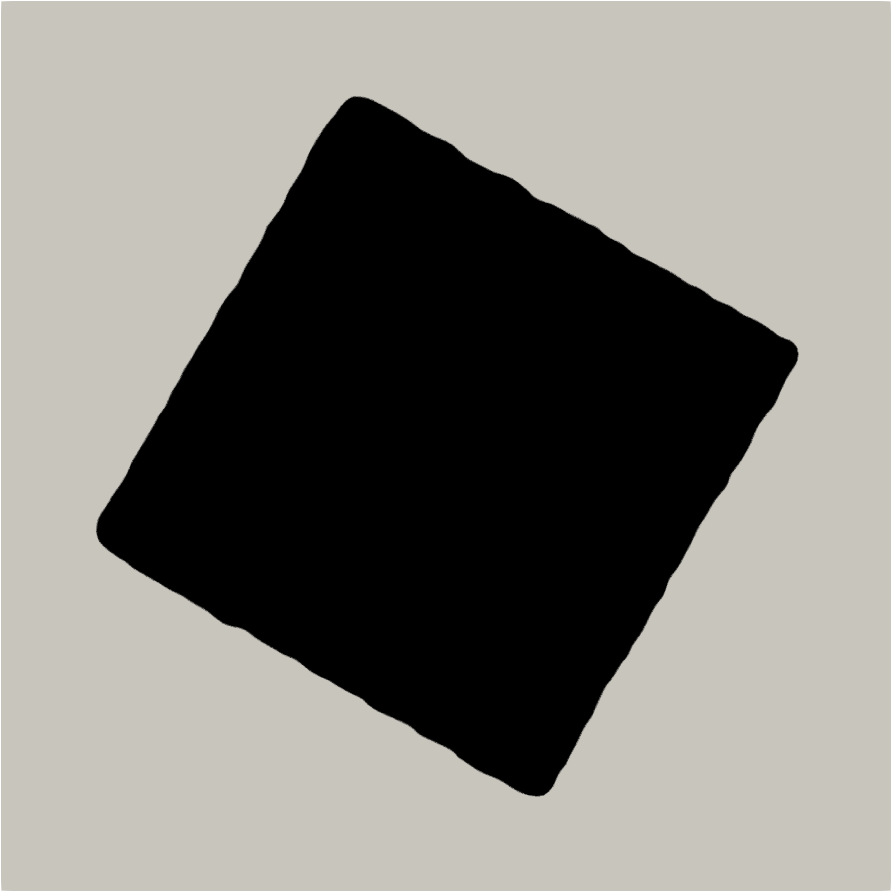}}	
	\subfigure[Final mesh]{\includegraphics[width=0.45\textwidth]{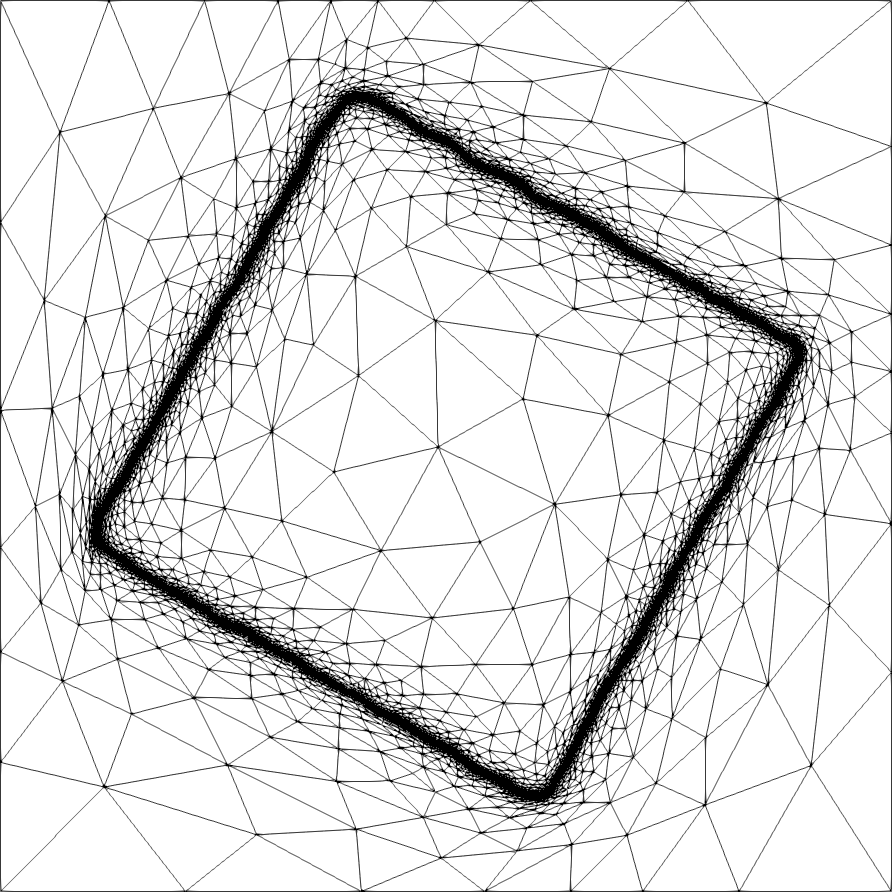}}
		
	\caption{Synthetic image with spatial information - Final (a) segmentation and (b) adapted mesh using the split-adapt Bregman algorithm with the Bayesian model.}
\label{fig:noiseMesh}
\end{figure}

\begin{figure}[!tb]
	\centering
	\subfigure[Split Bregman]{\includegraphics[width=0.45\textwidth]{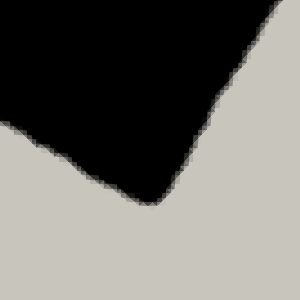}}
	\subfigure[Split-adapt Bregman]{\includegraphics[width=0.45\textwidth]{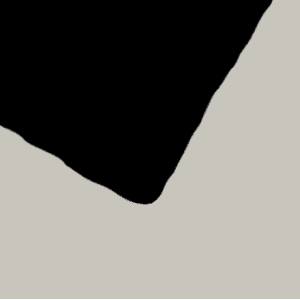}}	
		
	\caption{Synthetic image with spatial information - Detail of the final segmentation obtained using (a) the standard split Bregman and (b) the split-adapt Bregman algorithm with the Bayesian model.}
\label{fig:noiseDetail}
\end{figure}

Starting from the same initial guess and considering the same tolerance $\eta^\star$ for the stopping criterion of  the optimisation loop, the split-adapt Bregman method converges in 12 iterations with 4 mesh adaptation steps. Figure~\ref{fig:noiseMesh} displays the final segmentation and the corresponding adapted mesh consisting of 15,096 triangles, with $\hMin {=} 0.4 \times 10^{-1}$, $\hMax {=} 49.56$ and a maximum stretching factor $s_K {=} 847.3$. The split-adapt Bregman algorithm is thus capable of improving the description of the interface by avoiding mesh-dependent jagged effects, see figure~\ref{fig:noiseDetail}, while reducing the number of required mesh elements by almost $80 \%$. 

\subsection{Robustness to different noise sources}
\label{sc:Robustness}

In this section, the robustness of the split-adapt Bregman method is explored by considering the segmentation of three real images with inhomogeneous spatial information. The images are corrupted by introducing Gaussian, salt and pepper and speckle noise, respectively \hl{(see table~\ref{tab:ImageData}). Only the Bayesian segmentation model is considered, due to its superior performance with respect to the RSFE functional when complex patterns and spatial information are involved.  The input parameters for the split Bregman and split-adapt Bregman algorithms are detailed in~\ref{sc:appParam}.}
\begin{table}[!htb]
\begin{center}
\hl{
\begin{tabular}{c|c|c|c|c}
\hline
Image & Figure & Size (pixels) & Noise type & Noise level (variance) \\
  \hline
Tiger &   \ref{fig:tiger}  & $306 \times 212$ & Gaussian & $10^{-1}$ \\
   \hline
Boat & \ref{fig:boat} & $307 \times 259$ & Salt and pepper & $0.5  \times 10^{-1}$ \\
   \hline
Leaf &  \ref{fig:leaf} & $171 \times 135$  & Speckle & $0.1  \times 10^{-1}$ \\
   \hline
\end{tabular}
}
\end{center}
\caption{\hl{Real images - Details about the original image and the added noise.}}
\label{tab:ImageData}
\end{table}

\begin{figure}[!htb]
	\centering
	\subfigure[Noise-free image \label{fig:tigerImg}]{\includegraphics[width=0.32\textwidth]{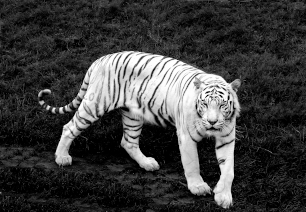}}
	\subfigure[Noisy image \label{fig:tigerNoiseImg}]{\includegraphics[width=0.32\textwidth]{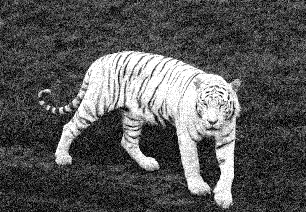}}	
	\subfigure[Ground truth \label{fig:tigerBayes}]{\includegraphics[width=0.32\textwidth]{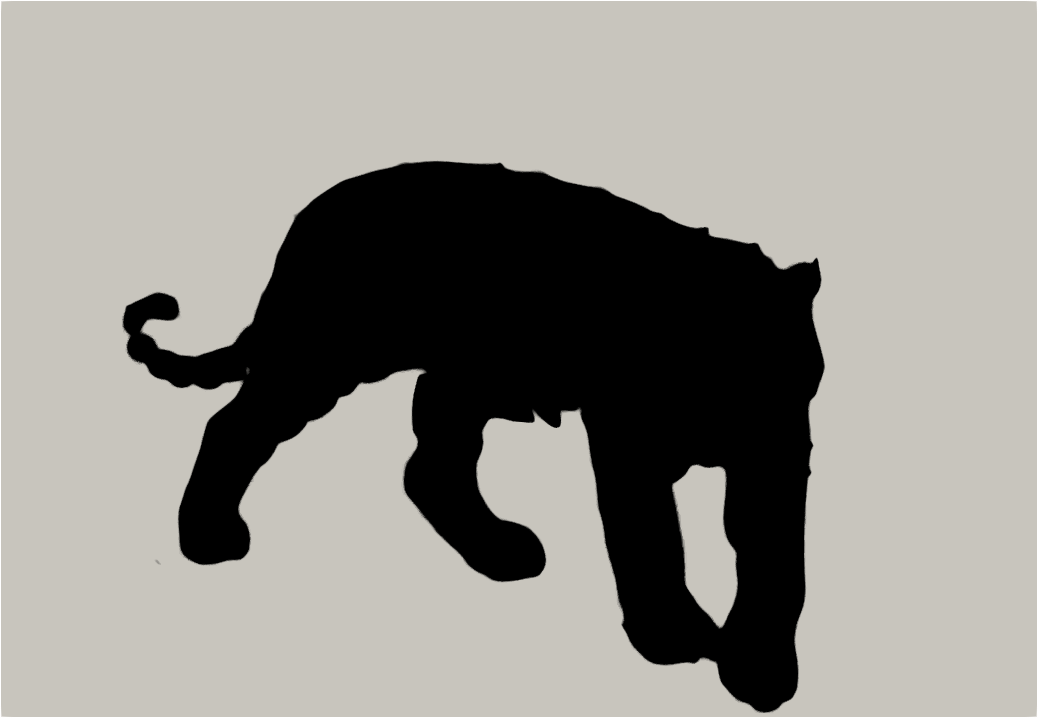}}
	
	\subfigure[Split Bregman \label{fig:tigerNoisy_Bayes}]{\includegraphics[width=0.32\textwidth]{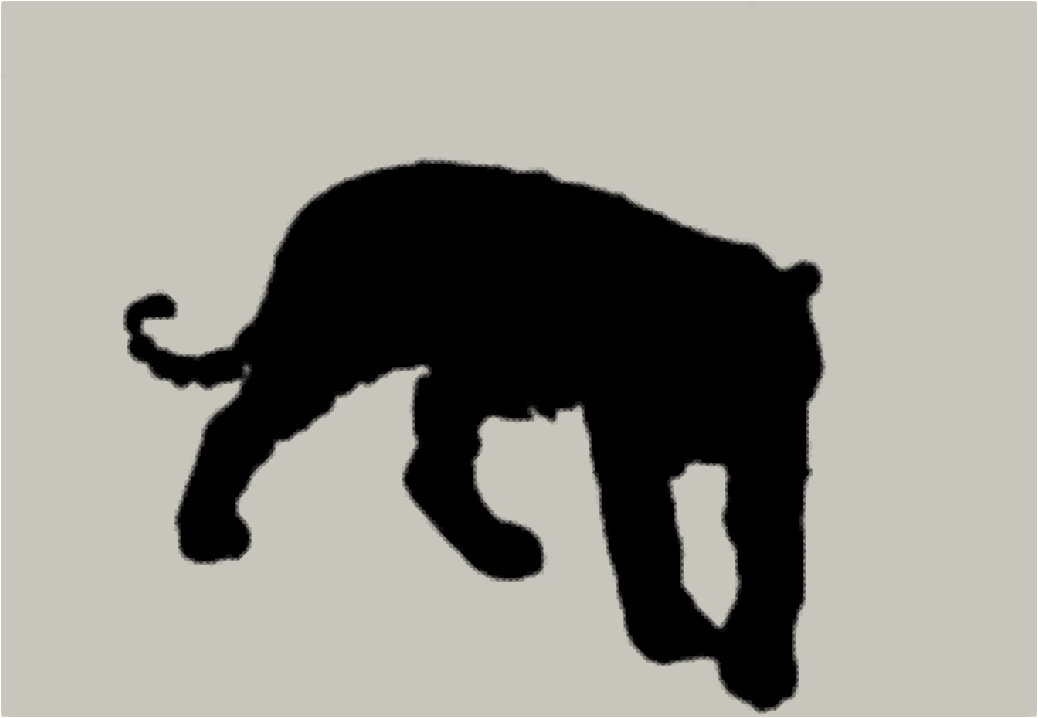}}
	\subfigure[Split-adapt Bregman \label{fig:tigerNoisy_Bayes_ad}]{\includegraphics[width=0.32\textwidth]{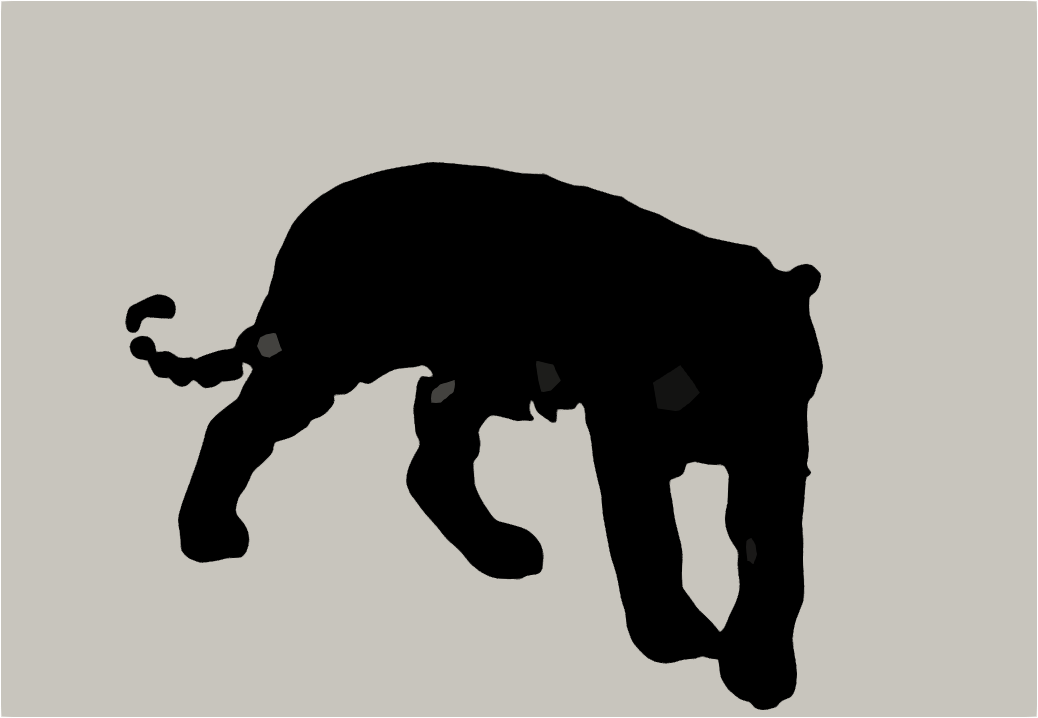}}	
	\subfigure[Final mesh \label{fig:tigerNoisy_mesh}]{\includegraphics[width=0.32\textwidth]{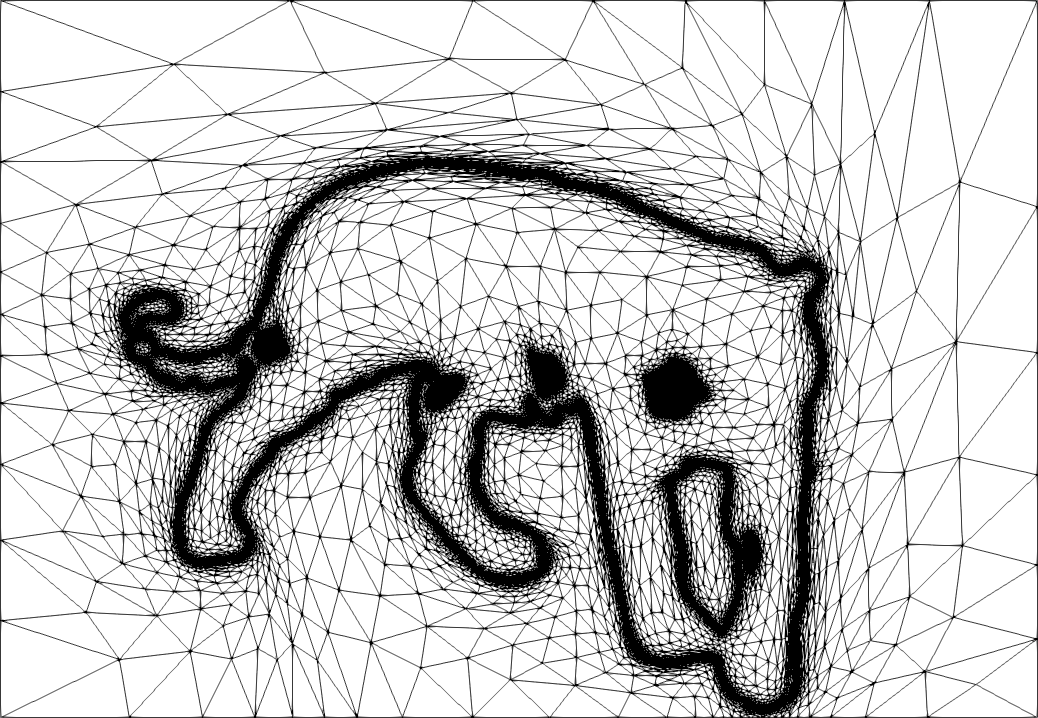}}
	
	\subfigure[Detail of \ref{fig:tigerNoisy_Bayes} \label{fig:tigerNoisy_Bayes_detail}]{\includegraphics[width=0.32\textwidth]{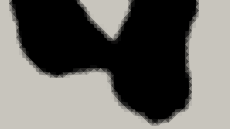}}
	\subfigure[Detail of \ref{fig:tigerNoisy_Bayes_ad} \label{fig:tigerNoisy_Bayes_ad_detail}]{\includegraphics[width=0.32\textwidth]{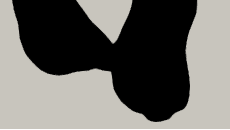}}	
	\subfigure[Detail of \ref{fig:tigerNoisy_mesh} \label{fig:tigerNoisy_mesh_detail}]{\includegraphics[width=0.32\textwidth]{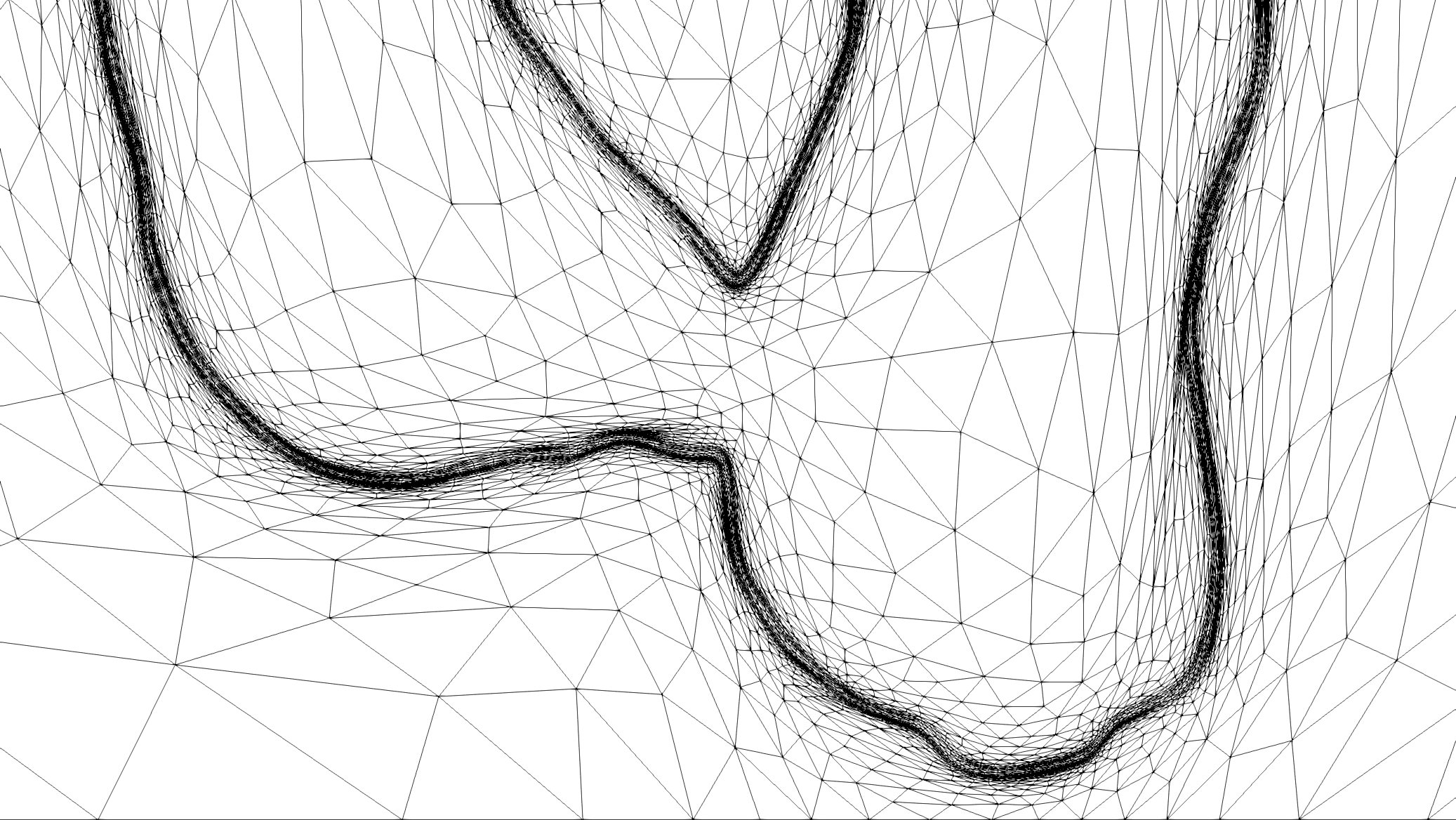}}
		
	\caption{Real image with Gaussian noise - (a) Noise-free image, (b) image corrupted by Gaussian noise with variance $0.1$ and (c) ground truth segmentation of the noise-free image provided by the split-adapt Bregman algorithm. Final segmented results of the noisy image computed using (d) the standard split Bregman and (e) the split-adapt Bregman algorithm. (f) Final adapted mesh and (g-i) details of the segmentation and the adapted mesh.}
\label{fig:tiger}
\end{figure}

\hl{As a first case study,  the image of the tiger in figure~\ref{fig:tigerImg} is corrupted by a Gaussian noise (see figure~\ref{fig:tigerNoiseImg}).  The details about the original noise-free image and the added noise are reported in table~\ref{tab:ImageData}.}
The ground truth segmentation is obtained using the split-adapt Bregman method to minimise the Bayesian functional on the noise-free image. \hl{Algorithm~\ref{alg:splitAdaptBregman} is executed starting from the data in table~\ref{tab:inputSetup}, which gathers the value of the tolerance $\eta^\star$ driving the segmentation procedure, together with the main features characterising the initial mesh (i.e., the number $\numel$ of mesh elements, the minimum element size $\hMin$ and the  maximum element size $\hMax$).
The split-adapt Bregman algorithm converges in 54 iterations with 18 mesh adaptation steps, providing the final segmentation displayed in figure~\ref{fig:tigerBayes} on an adapted mesh featuring 206,568 triangular elements, with $\hMin {=} 0.9 \times 10^{-2}$ and $\hMax {=} 76.77$, while the stretching factor achieves the maximum admissible value $s_K {=} 1,000$. }
\begin{table}[!htb]
\begin{center}
\hl{
\begin{tabular}{c|c|c|c|c}
\hline
Image & $\eta^\star$ & $\numel$ & $\hMin$ & $\hMax$ \\
  \hline
Tiger & $0.1 \times 10^{-2}$ & $128,710$ & $1$ & $1.4$ \\
   \hline
Boat & $0.5 \times 10^{-2}$ & $157,896$ & $1$ & $1.4$ \\
   \hline
Leaf & $0.1 \times 10^{-2}$ & $45,560$ & $1$ & $1.4$ \\
   \hline
\end{tabular}
}
\end{center}
\caption{\hl{Real images - Parameters for algorithm~\ref{alg:splitBregman} and for the initial configuration of algorithm~\ref{alg:splitAdaptBregman}.}}
\label{tab:inputSetup}
\end{table}

Both the split Bregman and the split-adapt Bregman methods provide accurate segmentations of the noisy image as confirmed by figures~\ref{fig:tigerNoisy_Bayes} and~\ref{fig:tigerNoisy_Bayes_ad}, respectively.  \hl{Computational details about the two algorithms are presented in table~\ref{tab:outputTime}, which collects the number $\nOpt$ of the optimisation iterations, the number $\nAdapt$ of mesh adaptation steps, the average execution time for one optimisation iteration ($\tOpt$) and for one mesh adaptation step ($\tAdapt$). Although the split-adapt Bregman algorithm requires more iterations than the split Bregman method to converge,  the average cost of each solve is reduced ($\SI{4.89}{s}$ versus $\SI{5.25}{s}$).  This is a direct consequence of the coarser meshes employed by the split-adapt Bregman strategy. A comparison between the uniform structured mesh used by the split Bregman algorithm and the anisotropically adapted one provided by the split-adapt Bregman method is presented in table~\ref{tab:outputMesh} in terms of the number $\numel$ of mesh elements,  the minimum ($\hMin$) and the maximum ($\hMax$) element size and the maximum element stretching factor $s_K$, which here achieves the admissible upper value limit $1,000$. The final mesh is displayed in figure~\ref{fig:tigerNoisy_mesh}.} Moreover, in figures~\ref{fig:tigerNoisy_Bayes_detail} and~\ref{fig:tigerNoisy_Bayes_ad_detail}, a detail of the two segmentations is shown, highlighting the superiority of the split-adapt Bregman method: this approach is able to improve the accuracy of  the interface representation, while reducing the number of the mesh elements by almost $40 \%$. Finally, figure~\ref{fig:tigerNoisy_mesh_detail} displays a detail of the adapted mesh, showing the benefit of employing elements with a large stretching factor to improve the smoothness and the accuracy of the description of the interface.
\begin{table}[!htb]
\begin{center}
\hl{
\begin{tabular}{c|c|c|c|c|c}
\hline
Image & Method &  $\nOpt$ & $\nAdapt$  & $\tOpt$ & $\tAdapt$ \\
\hline
Tiger & Split Bregman & $19$ & $-$ & $\SI{5.25}{s}$ & $-$ \\
\cline{2-6} & Split-adapt Bregman & $28$ & $9$ & $\SI{3.56}{s}$ & $\SI{1.33}{s}$ \\
\hline
Boat & Split Bregman & $9$ & $-$ & $\SI{6.64}{s}$ & $-$ \\
\cline{2-6} & Split-adapt Bregman & $10$ & $3$ & $\SI{2.02}{s}$ & $\SI{0.29}{s}$ \\
\hline
Leaf & Split Bregman & $20$ & $-$ & $\SI{1.91}{s}$ & $-$ \\
\cline{2-6} & Split-adapt Bregman & $22$ & $7$ & $\SI{2.27}{s}$ & $\SI{0.82}{s}$ \\
\hline
\end{tabular}
}
\end{center}
\caption{\hl{Real images - Comparison between algorithms~\ref{alg:splitBregman} and~\ref{alg:splitAdaptBregman} in terms of the computational cost.}}
\label{tab:outputTime}
\end{table}
\begin{table}[!htb]
\begin{center}
\hl{
\begin{tabular}{c|c|c|c|c|c}
\hline
Image & Method &  $\numel$ & $\hMin$  & $\hMax$ & $s_K$ \\
\hline
Tiger & Split Bregman & $128,710$ & $1$ & $1.4$ & $1$ \\
\cline{2-6} & Split-adapt Bregman & $77,197$ & $0.1 \times 10^{-1}$ & $67.45$ & $1,000$ \\
\hline
Boat & Split Bregman & $157,896$ & $1$ & $1.4$ & $1$ \\
\cline{2-6} & Split-adapt Bregman & $9,739$ & $0.56 \times 10^{-1}$ & $93.4$ & $482.93$ \\
\hline
Leaf & Split Bregman & $45,560$ & $1$ & $1.4$ & $1$ \\
\cline{2-6} & Split-adapt Bregman & $50,366$ & $0.15 \times 10^{-1}$ & $44.52$ & $1,000$ \\
\hline
\end{tabular}
}
\end{center}
\caption{\hl{Real images - Comparison between algorithms~\ref{alg:splitBregman} and~\ref{alg:splitAdaptBregman} in terms of the employed computational meshes.}}
\label{tab:outputMesh}
\end{table}

\begin{figure}[!htb]
	\centering
	\subfigure[Noise-free image~\ref{fig:tigerImg}]{\includegraphics[width=0.49\textwidth]{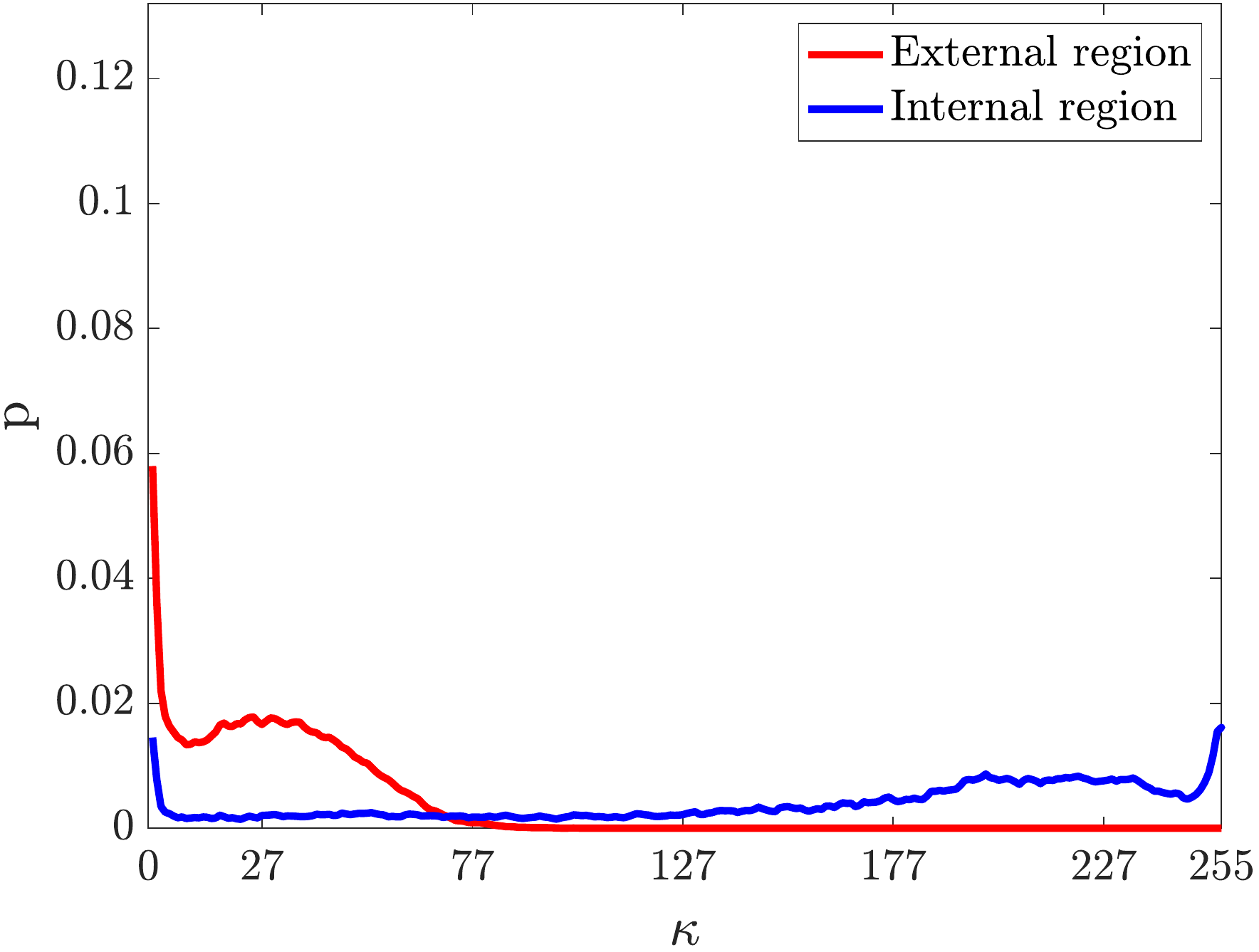}}
	\subfigure[Noisy image~\ref{fig:tigerNoiseImg}]{\includegraphics[width=0.49\textwidth]{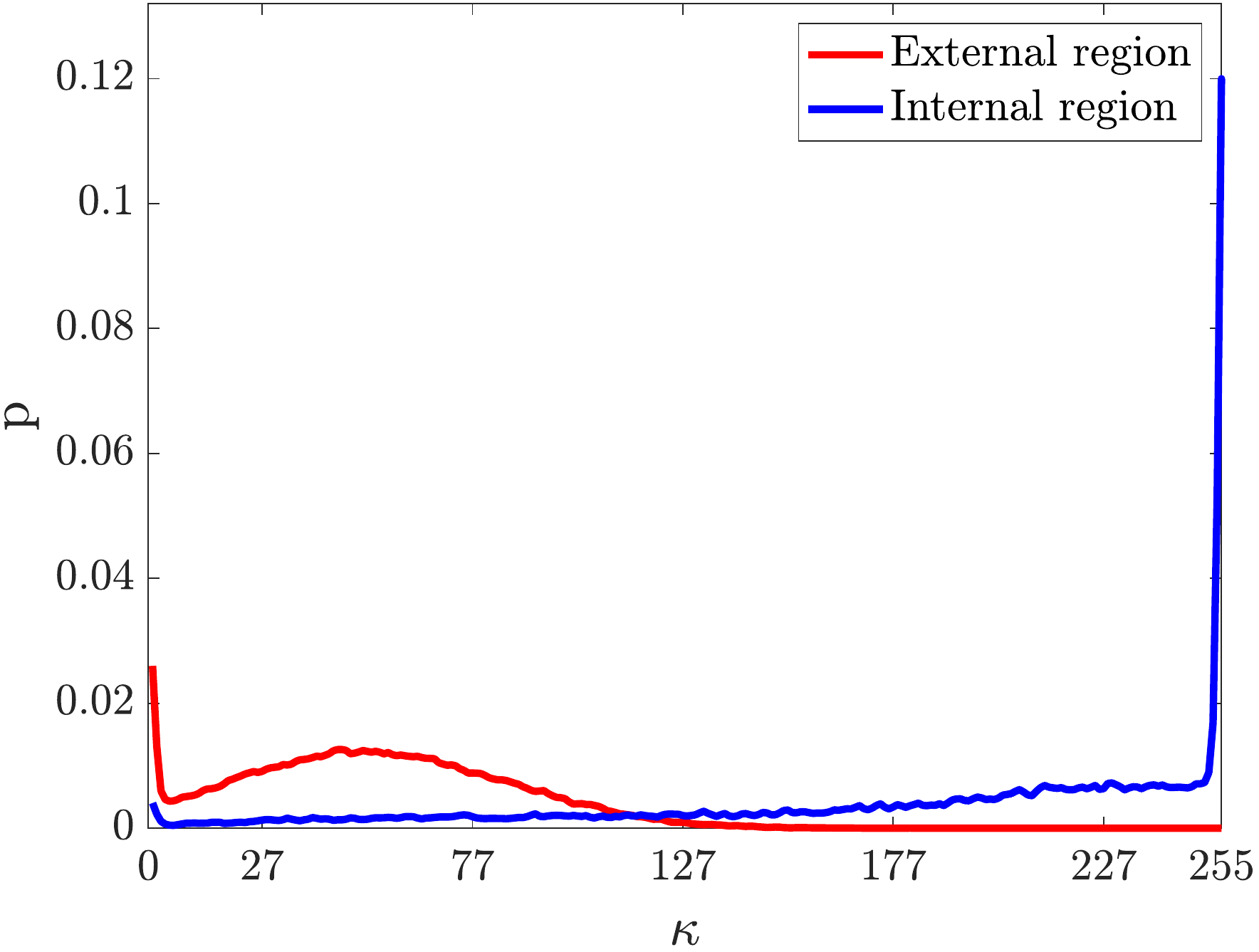}}	
		
	\caption{Real image with Gaussian noise - Estimated probability density functions of the pixel intensity $\kappa$ for the (a) noise-free and (b) noisy images.}
\label{fig:tigerDensity}
\end{figure}

It is worth noticing that the introduction of the noise into the image is responsible for a change in the spatial information as shown in figure~\ref{fig:tigerDensity}, where the estimated probability density distributions for the original and the noisy images are reported. Indeed, it is evident that the noise increases the variance of the background. More precisely, the values of the pixel intensity associated with the external region increases from the interval $[0,77]$ in the noise-free case to the interval $[0,127]$ in the noisy image, introducing additional difficulties for the segmentation algorithms. This is particularly evident in localised regions where the difference between the background and the foreground is less clear (e.g., the tail) and the corresponding segmentation suffers from a slight loss of accuracy.

\begin{figure}[!htb]
	\centering
	\subfigure[Noise-free image \label{fig:boatImg}]{\includegraphics[width=0.32\textwidth]{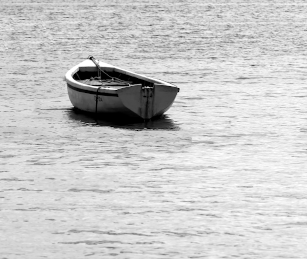}}
	\subfigure[Noisy image \label{fig:boatNoiseImg}]{\includegraphics[width=0.32\textwidth]{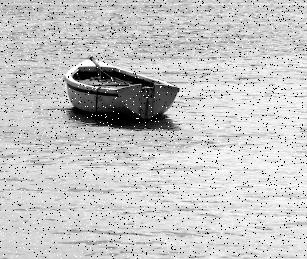}}	
	\subfigure[Ground truth \label{fig:boatBayes}]{\includegraphics[width=0.32\textwidth]{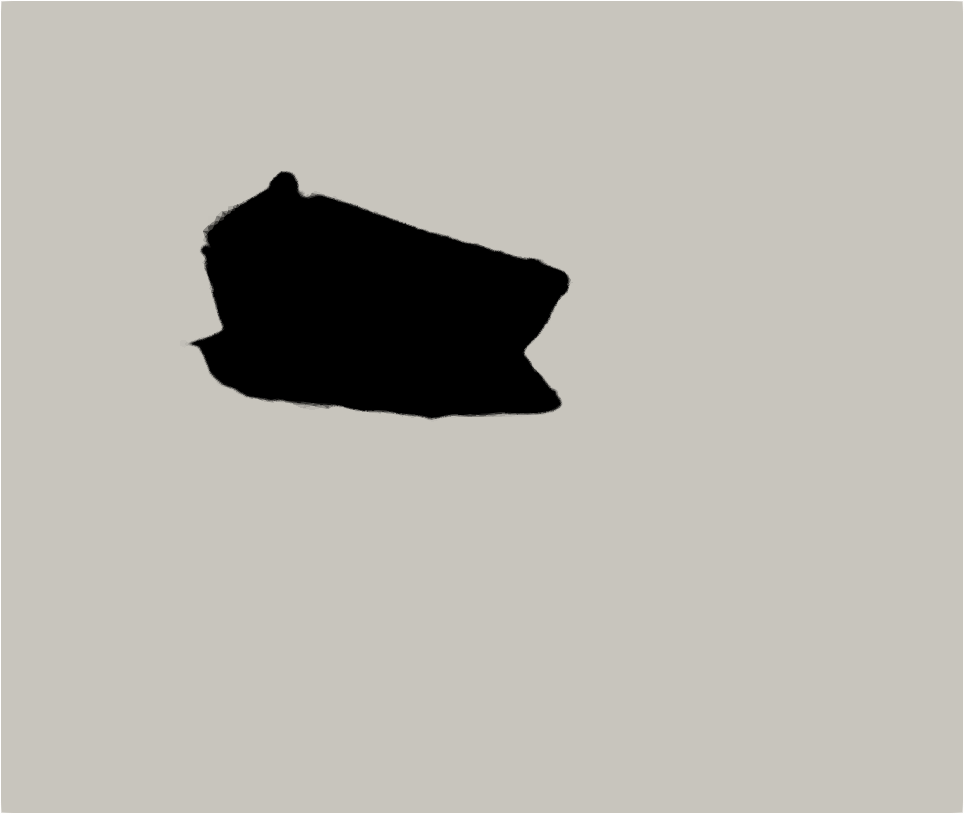}}
	
	\subfigure[Split Bregman \label{fig:boatNoisy_Bayes}]{\includegraphics[width=0.32\textwidth]{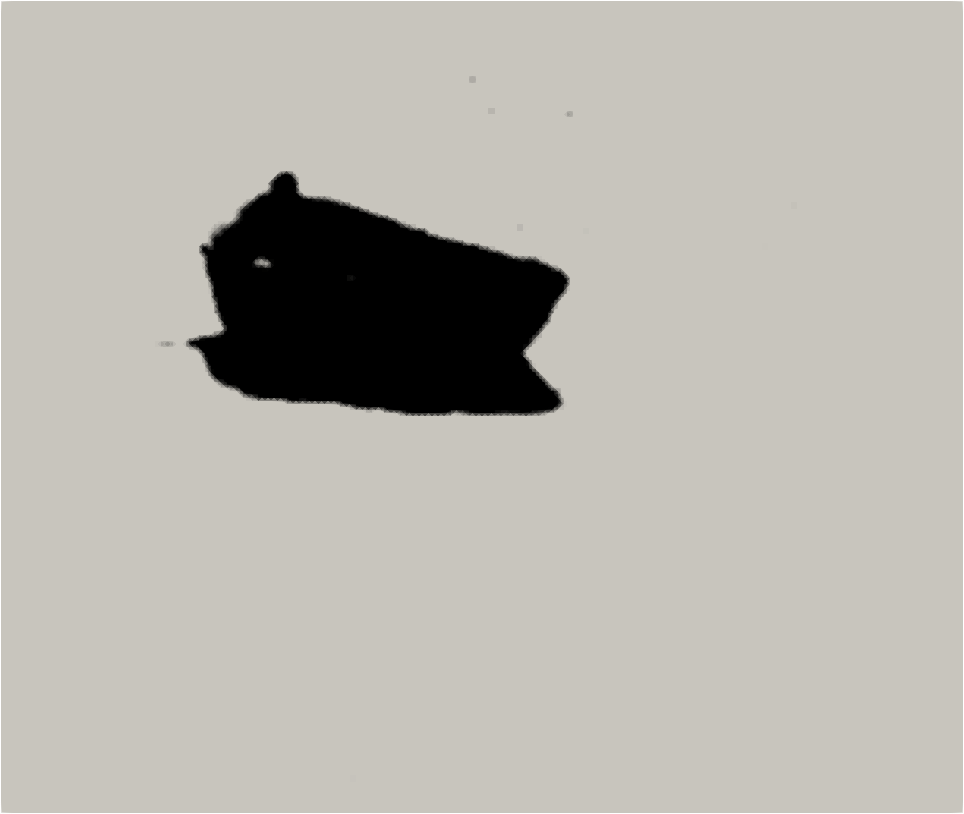}}
	\subfigure[Split-adapt Bregman \label{fig:boatNoisy_Bayes_ad}]{\includegraphics[width=0.32\textwidth]{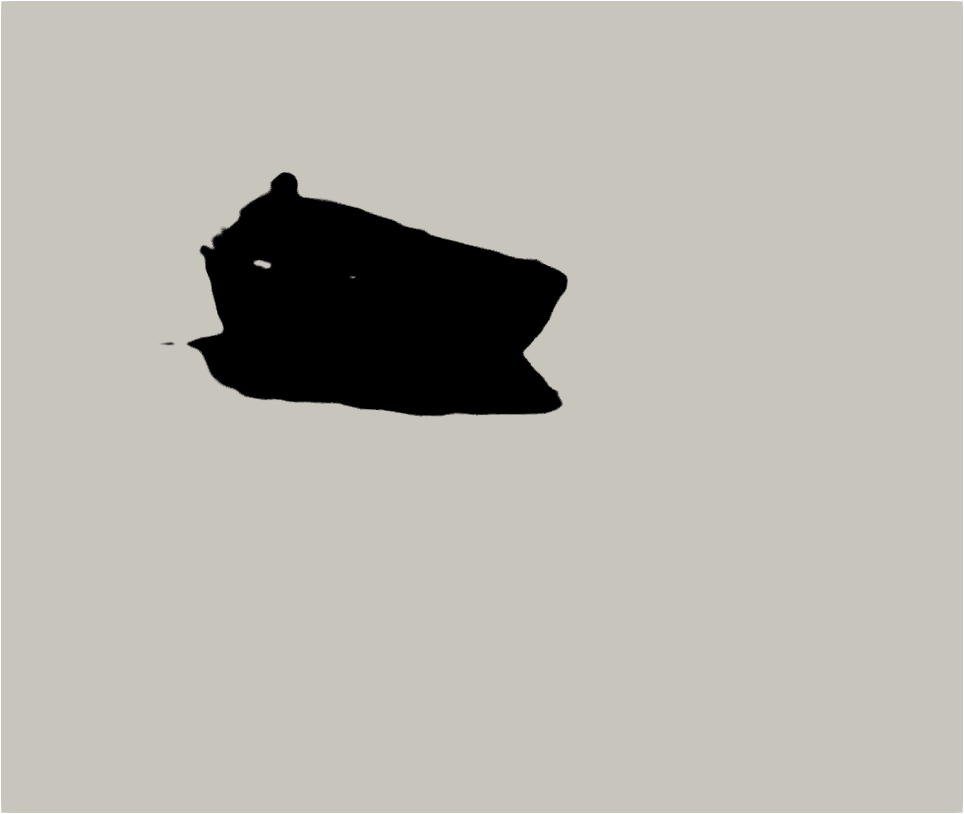}}	
	\subfigure[Final mesh \label{fig:boatNoisy_mesh}]{\includegraphics[width=0.32\textwidth]{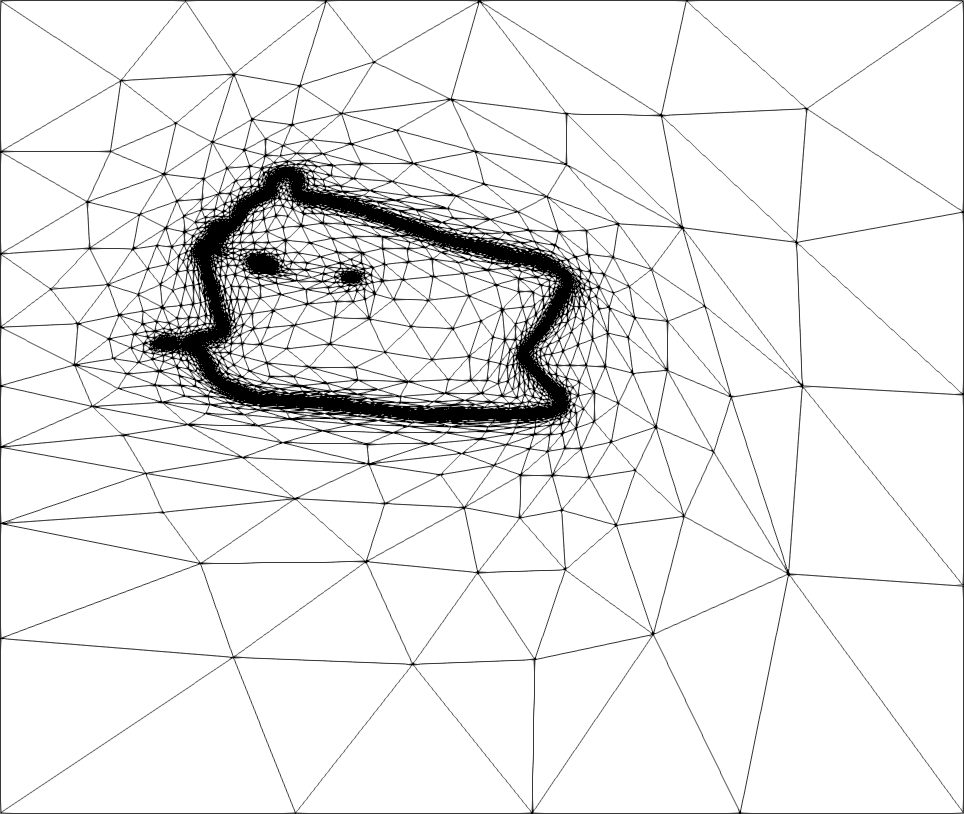}}
	
	\subfigure[Detail of \ref{fig:boatNoisy_Bayes} \label{fig:boatNoisy_Bayes_detail}]{\includegraphics[width=0.32\textwidth]{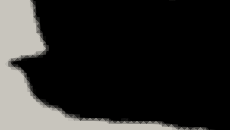}}
	\subfigure[Detail of \ref{fig:boatNoisy_Bayes_ad} \label{fig:boatNoisy_Bayes_ad_detail}]{\includegraphics[width=0.32\textwidth]{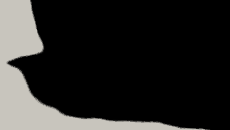}}	
	\subfigure[Detail of \ref{fig:boatNoisy_mesh} \label{fig:boatNoisy_mesh_detail}]{\includegraphics[width=0.32\textwidth]{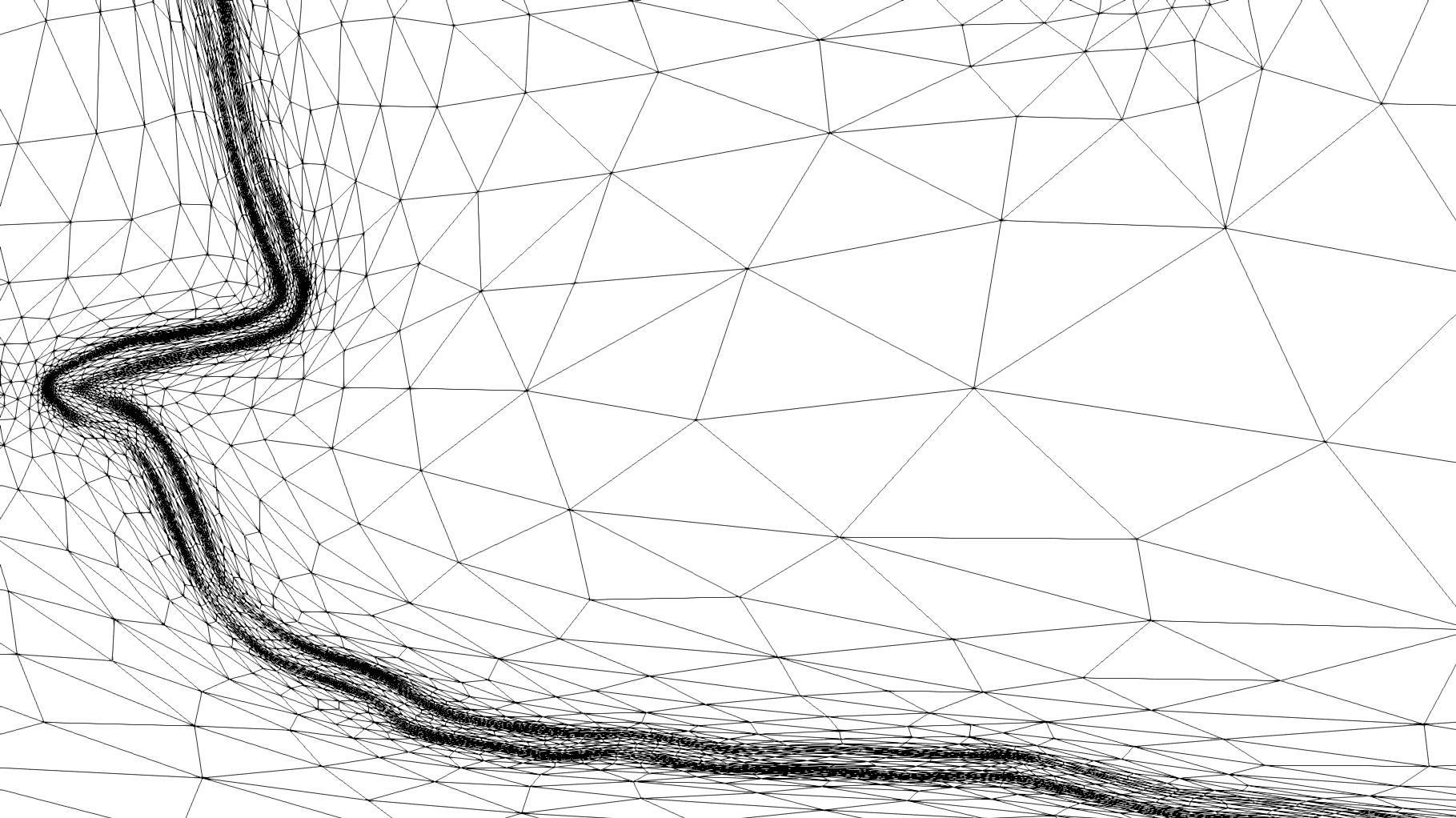}}	
	\caption{Real image with salt and pepper noise - (a) Noise-free image, (b) image corrupted by salt and pepper noise with variance $0.5 \times 10^{-1}$ and (c) ground truth segmentation of the noise-free image provided by the split-adapt Bregman algorithm. Final segmented results of the noisy image computed using (d) the standard split Bregman and (e) the split-adapt Bregman algorithm. (f) Final adapted mesh and (g-i) details of the segmentation and the adapted mesh.}
\label{fig:boat}
\end{figure}

\hl{The second example features the image of the boat in figure~\ref{fig:boatImg}, corrupted by a salt and pepper noise as displayed in figure~\ref{fig:boatNoiseImg}.  Details on the noise level are presented in table~\ref{tab:ImageData}. Figure~\ref{fig:boatBayes} displays the ground truth segmentation provided by the split-adapt Bregman method applied to the noise-free image,  using the initial mesh and the tolerance $\eta^\star$ detailed in table~\ref{tab:inputSetup}.} The final mesh, obtained after 6 iterations of the split-adapt Bregman algorithm with 2 mesh adaptation steps, features 39,064 triangles, with $\hMin {=} 0.85 \times 10^{-1}$, $\hMax {=} 4.72$ and a maximum elemental stretching factor $s_K {=} 557.29$.
 
Both the split Bregman and the split-adapt Bregman algorithms confirm their capability of exploiting the spatial information to accurately segment noisy images, even in the presence of salt and pepper noise.  \hl{More precisely, figure~\ref{fig:boatNoisy_Bayes} shows the output of the split Bregman algorithm on a uniform structured mesh, whereas the result of the split-adapt Bregman approach, displayed in figure~\ref{fig:boatNoisy_Bayes_ad}, is obtained using the final adapted mesh in figure~\ref{fig:boatNoisy_mesh}.  Quantitative information about the output meshes is reported in table~\ref{tab:outputMesh},  whereas table~\ref{tab:outputTime} highlights the computational details of the two algorithms. As observed for the previous case, the split-adapt Bregman strategy converges in more iterations than the split Bregman method. However, the extra iterations introduce an almost negligible computational burden given the extremely coarse mesh designed by the adaptation tool.  Indeed, the average cost of each optimisation iteration drops from $\SI{6.64}{s}$ for the split Bregman method to $\SI{2.02}{s}$ for the split-adapt Bregman algorithm. The extra cost of each adaptation step is very limited since both the gradient reconstruction and the metric derivation are characterised by explicit formulae so that the generation of the new mesh lasts on average only $\SI{0.29}{s}$}. In summary, excellent results are obtained using a computational mesh almost $94 \%$ coarser than the uniform grid employed by the standard split Bregman method. Moreover, from the comparison of figures~\ref{fig:boatNoisy_Bayes_detail} and~\ref{fig:boatNoisy_Bayes_ad_detail}, it is evident the improvement led by mesh adaptation, with a limited computational effort, in terms of the description of the boundaries of the segmented region. A detail of the final adapted mesh is provided in figure~\ref{fig:boatNoisy_mesh_detail}, showing the local refinement of the elements and their anisotropy along the interface between the object and the background. 
 
\begin{figure}[!htb]
	\centering
	\subfigure[Noise-free image \label{fig:leafImg}]{\includegraphics[width=0.32\textwidth]{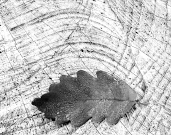}}
	\subfigure[Noisy image \label{fig:leafNoiseImg}]{\includegraphics[width=0.32\textwidth]{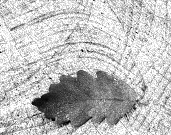}}	
	\subfigure[Ground truth \label{fig:leafBayes}]{\includegraphics[width=0.32\textwidth]{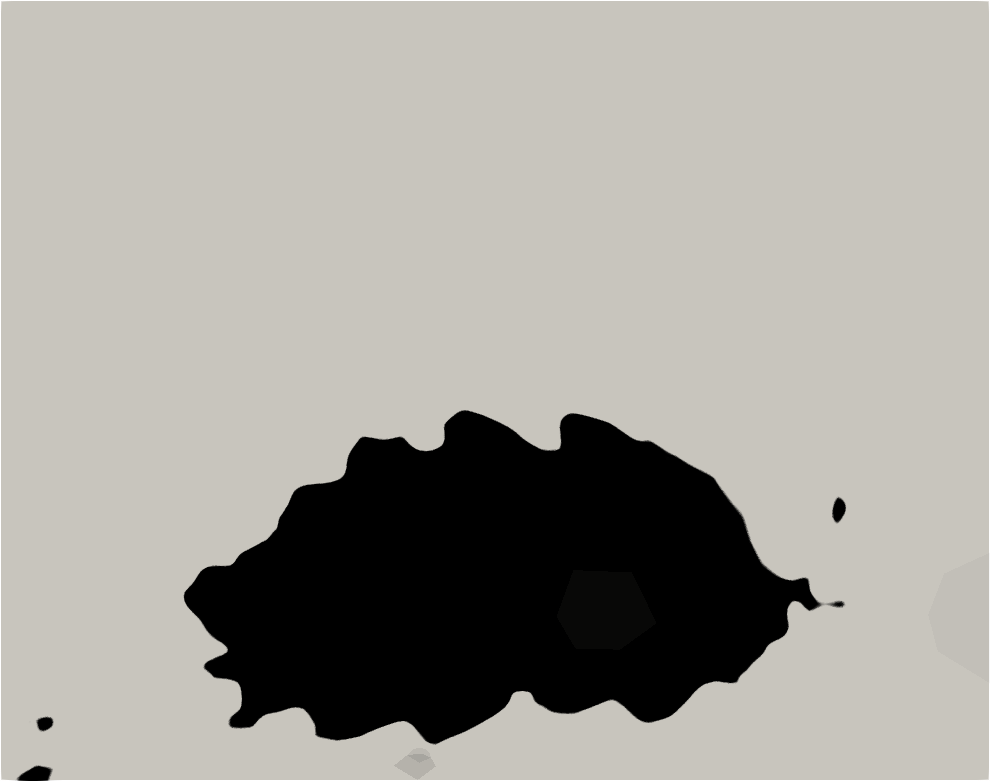}}
	
	\subfigure[Split Bregman \label{fig:leafNoisy_Bayes}]{\includegraphics[width=0.32\textwidth]{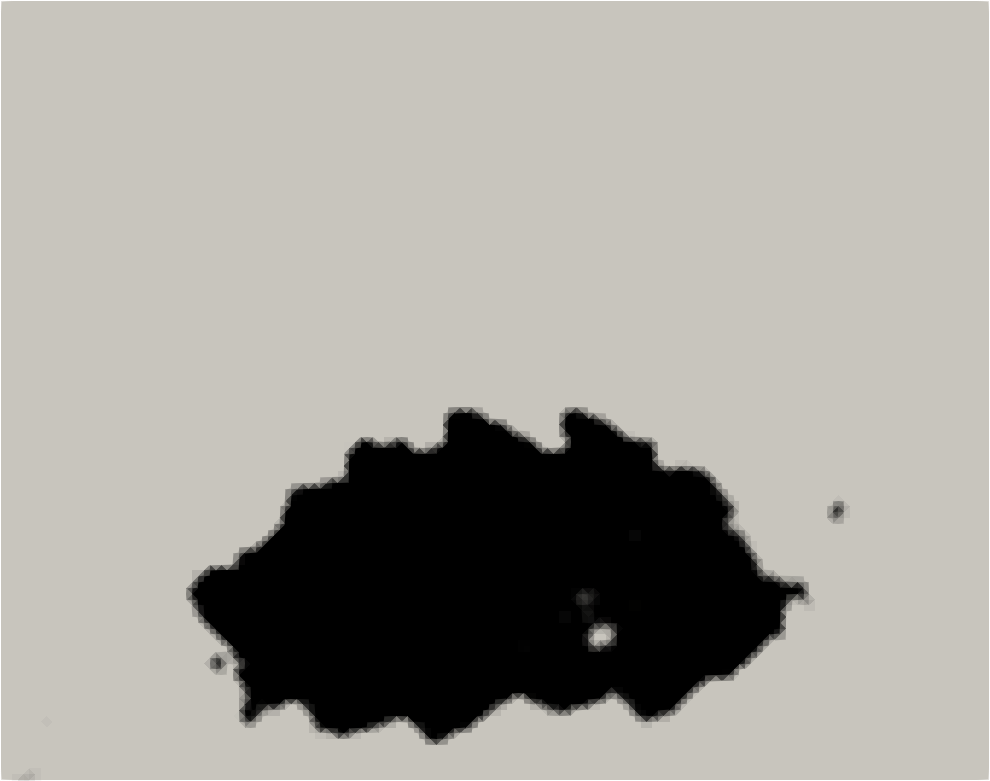}}
	\subfigure[Split-adapt Bregman \label{fig:leafNoisy_Bayes_ad}]{\includegraphics[width=0.32\textwidth]{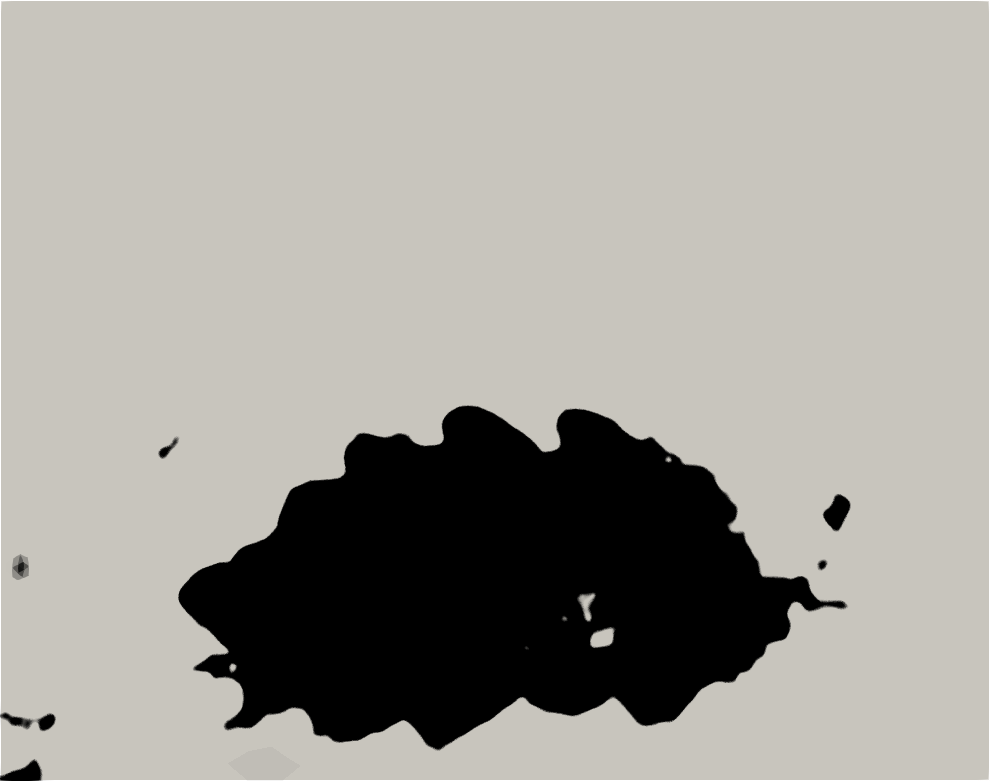}}	
	\subfigure[Final mesh \label{fig:leafNoisy_mesh}]{\includegraphics[width=0.32\textwidth]{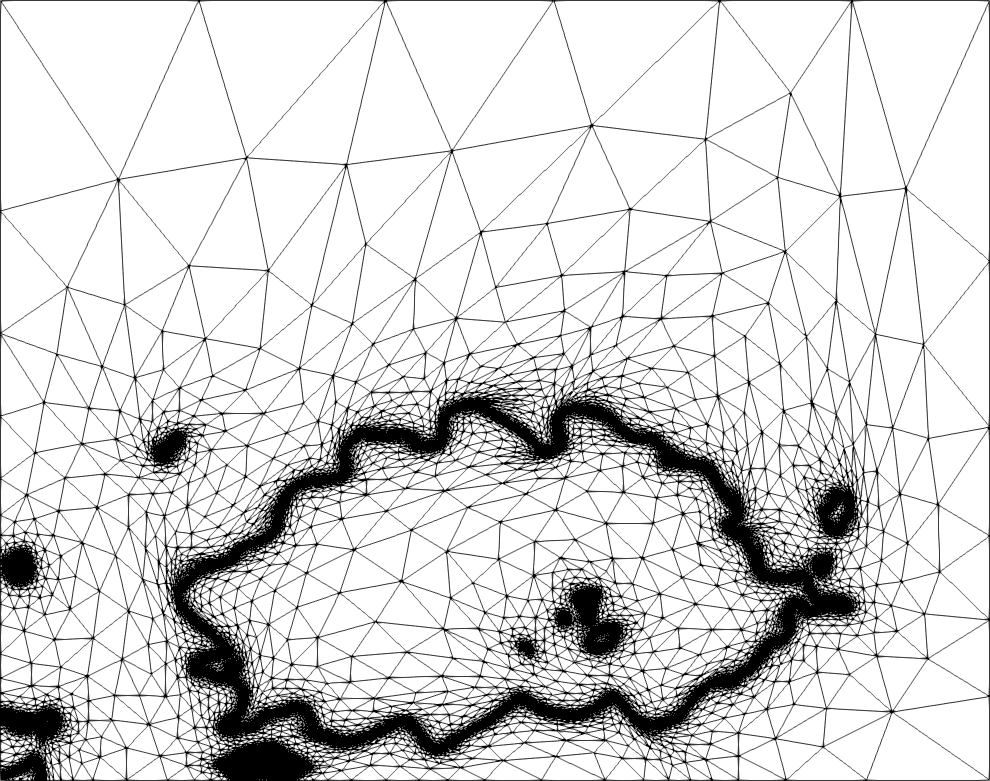}}
	
	\subfigure[Detail of \ref{fig:leafNoisy_Bayes} \label{fig:leafNoisy_Bayes_detail}]{\includegraphics[width=0.32\textwidth]{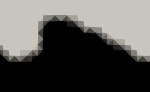}}
	\subfigure[Detail of \ref{fig:leafNoisy_Bayes_ad} \label{fig:leafNoisy_Bayes_ad_detail}]{\includegraphics[width=0.32\textwidth]{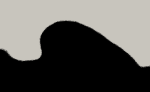}}	
	\subfigure[Detail of \ref{fig:leafNoisy_mesh} \label{fig:leafNoisy_mesh_detail}]{\includegraphics[width=0.32\textwidth]{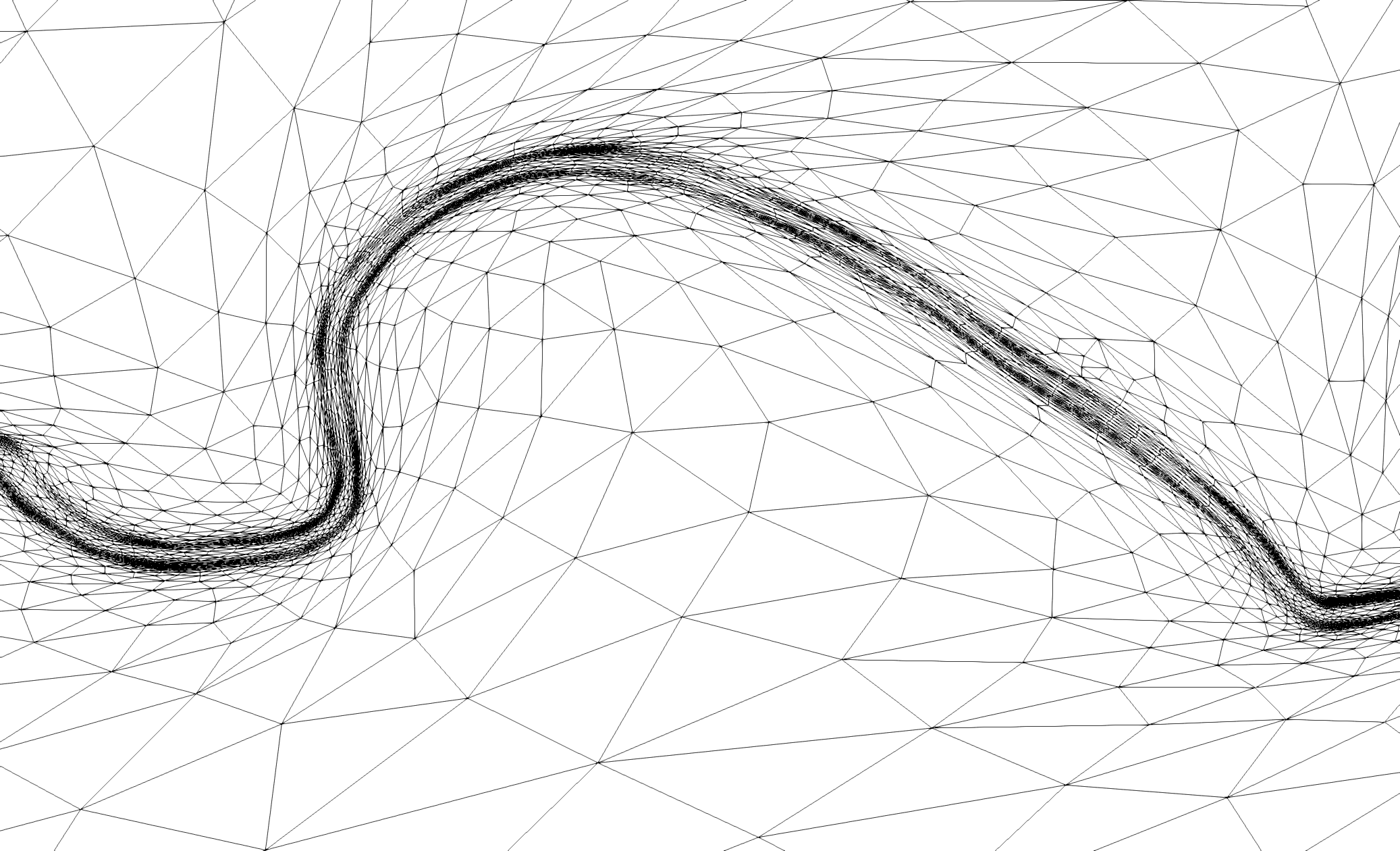}}
		
	\caption{Real image with speckle noise - (a) Noise-free image, (b) image corrupted by salt and pepper noise with variance $0.1 \times 10^{-1}$ and (c) ground truth segmentation of the noise-free image provided by the split-adapt Bregman algorithm. Final segmented results of the noisy image computed using (d) the standard split Bregman and (e) the split-adapt Bregman algorithm. (f) Final adapted mesh and (g-i) details of the segmentation and the adapted mesh.}
\label{fig:leaf}
\end{figure}

The last example of this section represents a challenge for segmentation algorithms, the goal being the treatment of an image corrupted by speckle noise.  \hl{The image of the leaf in figure~\ref{fig:leafImg} is corrupted  by a speckle noise, leading to the noisy version in figure~\ref{fig:leafNoiseImg} (see table~\ref{tab:ImageData} for additional details). Using the setup in table~\ref{tab:inputSetup} for the segmentation of the noise-free image, the split-adapt Bregman method yields the ground truth solution shown in figure~\ref{fig:leafBayes}. } The final adapted mesh features 40,491 elements with $\hMin {=} 0.88 \times 10^{-2}$ and $\hMax {=} 48.64$, while the stretching factor attains the maximum admissible value $s_K {=} 1,000$. The case under analysis is particularly challenging for segmentation algorithms given the richness of spatial information present in both the background and the foregound of the image. Indeed, the split-adapt Bregman method requires 84 iterations with 28 mesh adaptation steps to converge to the ground truth segmentation of the noise-free image. In addition, some numerical artifacts are detected in the form of small segmented regions near the boundary, as displayed in figure~\ref{fig:leafBayes}.
 
The image corrupted by speckle noise is segmented using the split Bregman and split-adapt Bregman algorithms.  \hl{The former approach is able to identify the region to be segmented, with minor numerical artifacts, which are localised in highly inhomogeneous areas, see figure~\ref{fig:leafNoisy_Bayes}. The split-adapt Bregman algorithm improves the overall representation of the boundary of the segmentation as highlighted by figure~\ref{fig:leafNoisy_Bayes_ad} and the detail in figure~\ref{fig:leafNoisy_Bayes_ad_detail}, although additional small numerical artifacts appear.  The configuration here considered is especially challenging because of the multiplicative nature of the speckle noise whose value depends locally on the pixel intensity.  This leads to artificial localised large gradients of the pixel intensity, which are responsible for triggering the mesh adaptation procedure based on the control of the gradient of the solution, with a consequent local refinement in the neighbourhood of the artifacts (see figure~\ref{fig:leafNoisy_mesh}). Moreover, the selected minimum element size also plays a role in the presence of such artifacts. As shown in table~\ref{tab:outputMesh},  the minimum element size is $\hMin {=} 0.15 \times 10^{-1}$ for the adapted mesh, whereas $\hMin {=} 1$ in the uniform structured grid employed by the split Bregman method. The larger value of $\hMin$ characterising the standard algorithm is responsible for the loss of all the details below a certain dimension, thus indirectly yielding a local regularisation of the segmented region.  Finally,  the overall quality of the solution in figure~\ref{fig:leafNoisy_Bayes_ad} {for the noisy image is fully comparable with the one of the ground truth segmentation of the noise-free image in figure~\ref{fig:leafBayes}, confirming the robustness of the split-adapt Bregman method even in the presence of speckle noise. On the contrary, the segmentation in figure~\ref{fig:leafNoisy_Bayes} (see also a detail in figure~\ref{fig:leafNoisy_Bayes_detail}) is definitely more blurred than the corresponding output of algorithm~\ref{alg:splitAdaptBregman}.
} 

The number of optimisation iterations required by algorithms~\ref{alg:splitBregman} and~\ref{alg:splitAdaptBregman} to converge is provided in table~\ref{tab:outputTime}, whereas the details on the employed computational meshes are furnished in table~\ref{tab:outputMesh}, where the maximum admissible stretching factor $s_K {=} 1,000$ is achieved by the split-adapt Bregman method.} Although in this case the cardinality of the adapted mesh is slightly superior to the one of the initial structured grid, the split-adapt Bregman algorithm still allows to improve the quality of the interface as it is clearly visible in figures~\ref{fig:leafNoisy_Bayes_detail} and~\ref{fig:leafNoisy_Bayes_ad_detail}. Indeed, the quality of the segmented boundary greatly improves owing to the local refinement and the anisotropy of the mesh elements, see figure~\ref{fig:leafNoisy_mesh_detail} for a detail.
 
The above simulations highlight the capability of the split-adapt Bregman method applied to the Bayesian segmentation model of exploiting the spatial information in the images under analysis, even when they are corrupted by different types of noise. In addition, the split-adapt Bregman method outperforms the standard split Bregman algorithm by significantly improving the quality of the boundaries of the segmented region while reducing the number of required degrees of freedom of the problem thanks to the proper alignment, sizing and shaping of the mesh elements.

\subsection{Application to medical images}
\label{sc:Ultrasound}

\begin{figure}[!htb]
	\centering
	\subfigure[Image \label{fig:ultrasoundImg}]{\includegraphics[width=0.45\textwidth]{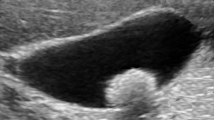}}
	\subfigure[Split Bregman \label{fig:ultrasound_Bayes}]{\includegraphics[width=0.45\textwidth]{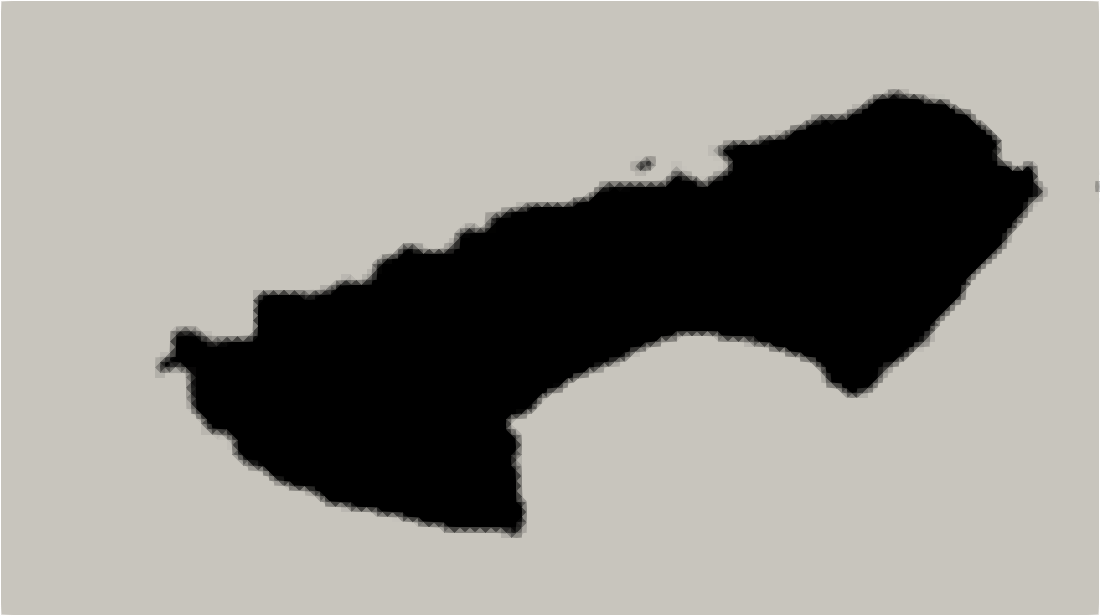}}	
	
	\subfigure[Split-adapt Bregman \label{fig:ultrasound_Bayes_ad}]{\includegraphics[width=0.45\textwidth]{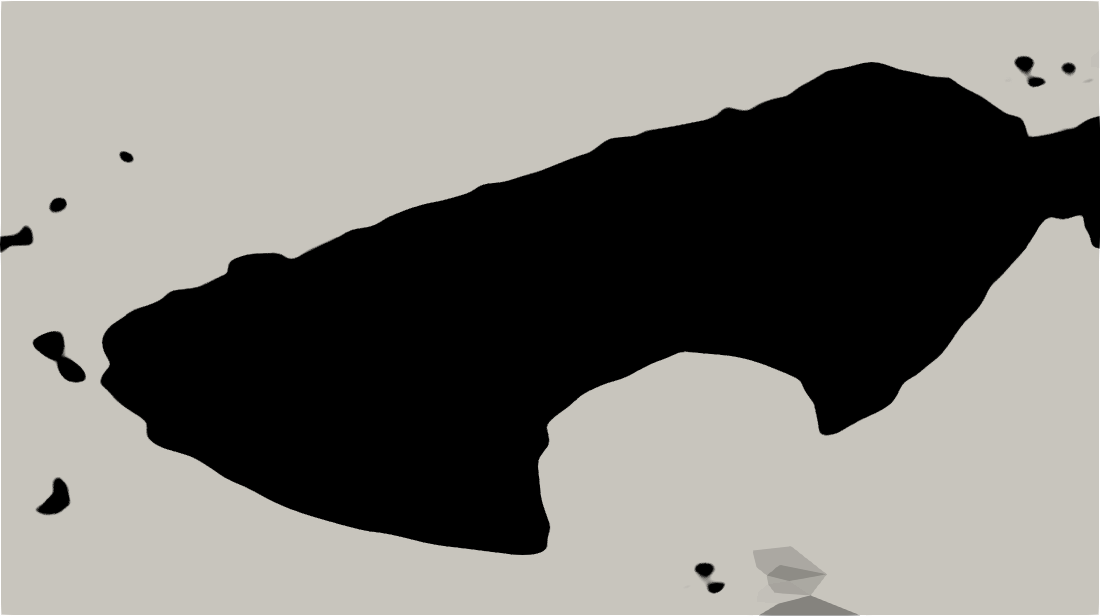}}
	\subfigure[Final mesh \label{fig:ultrasound_mesh}]{\includegraphics[width=0.45\textwidth]{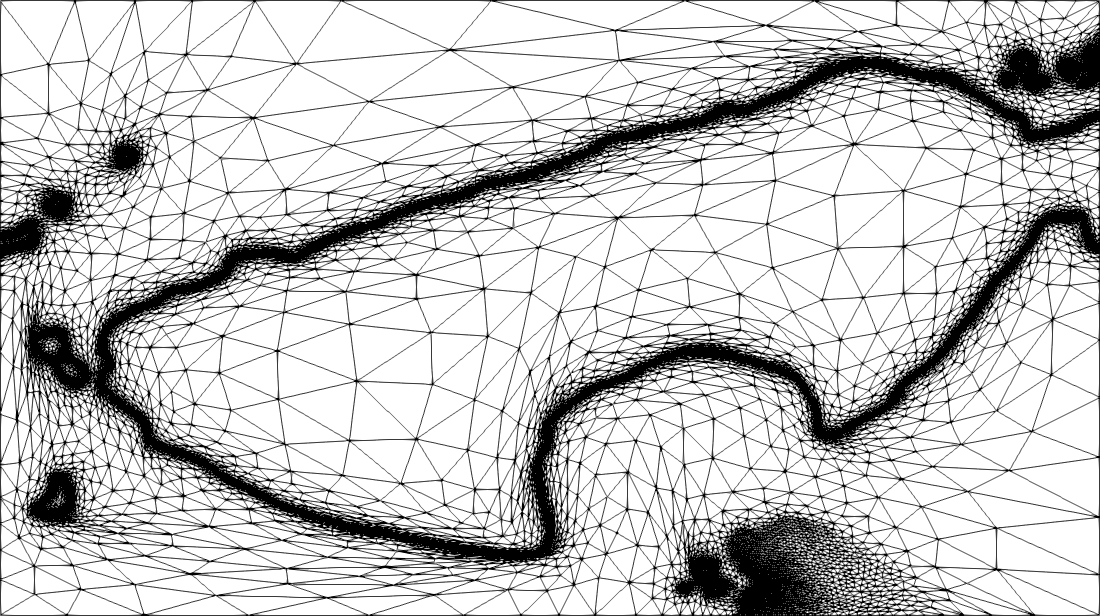}}	
	
	\subfigure[Detail of \ref{fig:ultrasound_Bayes_ad} \label{fig:ultrasound_Bayes_ad_detail}]{\includegraphics[width=0.45\textwidth]{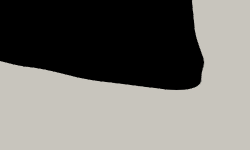}}
	\subfigure[Detail of \ref{fig:ultrasound_mesh} \label{fig:ultrasound_mesh_detail}]{\includegraphics[width=0.45\textwidth]{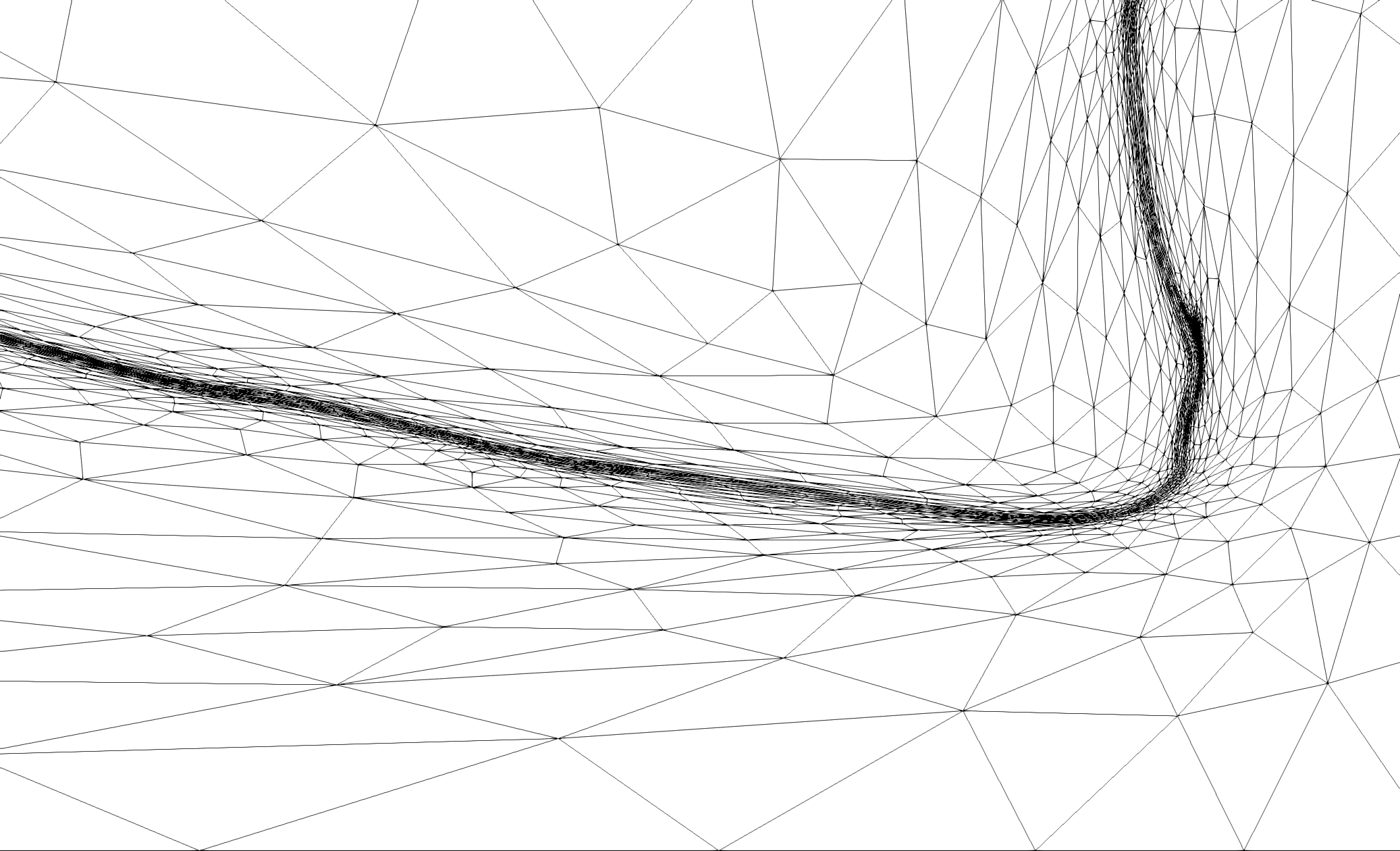}}	
		
	\caption{Real ultrasound image - (a) Image. Final segmented results computed using (b) the split Bregman and (c) the split-adapt Bregman algorithm. (d) Final adapted mesh and (e-f) details of the segmentation and the adapted mesh.}
\label{fig:ultrasound}
\end{figure}

The last example considers the segmentation of a $214 \times  120$ pixels medical image obtained via the ultrasound of a gallbladder, see figure~\ref{fig:ultrasoundImg}. Notice that the image is extremely noisy, thus constituting an actual challenge for segmentation algorithms.

The segmentation is performed using the split Bregman and split-adapt Bregman algorithms with the tolerance $\eta^\star {=} 0.1 \times 10^{-2}$ to stop the optimisation procedure.  \hl{A detailed description of the input parameters for the method is provided in~\ref{sc:appParam}.} The split Bregman method is executed on a uniform structured grid with 50,694 triangles such that $\hMin {=} 1$ and $\hMax {=} 1.4$ and it converges in 17 iterations. The split-adapt Bregman algorithm performs 94 optimisation iterations with 31 mesh adaptation steps to converge on the adapted mesh displayed in figure~\ref{fig:ultrasound_mesh}.  The final mesh consists of 89,845 triangular elements, with $\hMin {=} 0.95 \times 10^{-2}$ and $\hMax {=} 30.61$ and with the stretching factor achieving the maximum admissible value $s_K {=} 1,000$. Despite the same tolerance $\eta^\star$ is employed for both approaches, the split-adapt Bregman algorithm outperforms the standard method providing a more complete description of the region to be segmented, while capturing also finer features of the image. The resulting segmentations are shown in figure~\ref{fig:ultrasound_Bayes} and~\ref{fig:ultrasound_Bayes_ad} for comparison. In addition, consistently with previous simulations, the split-adapt Bregman method also improves the quality of the interface as confirmed by the details in figures~\ref{fig:ultrasound_Bayes_ad_detail} and~\ref{fig:ultrasound_mesh_detail}.

\begin{SCfigure}[][!htb]
	\includegraphics[width=0.6\textwidth]{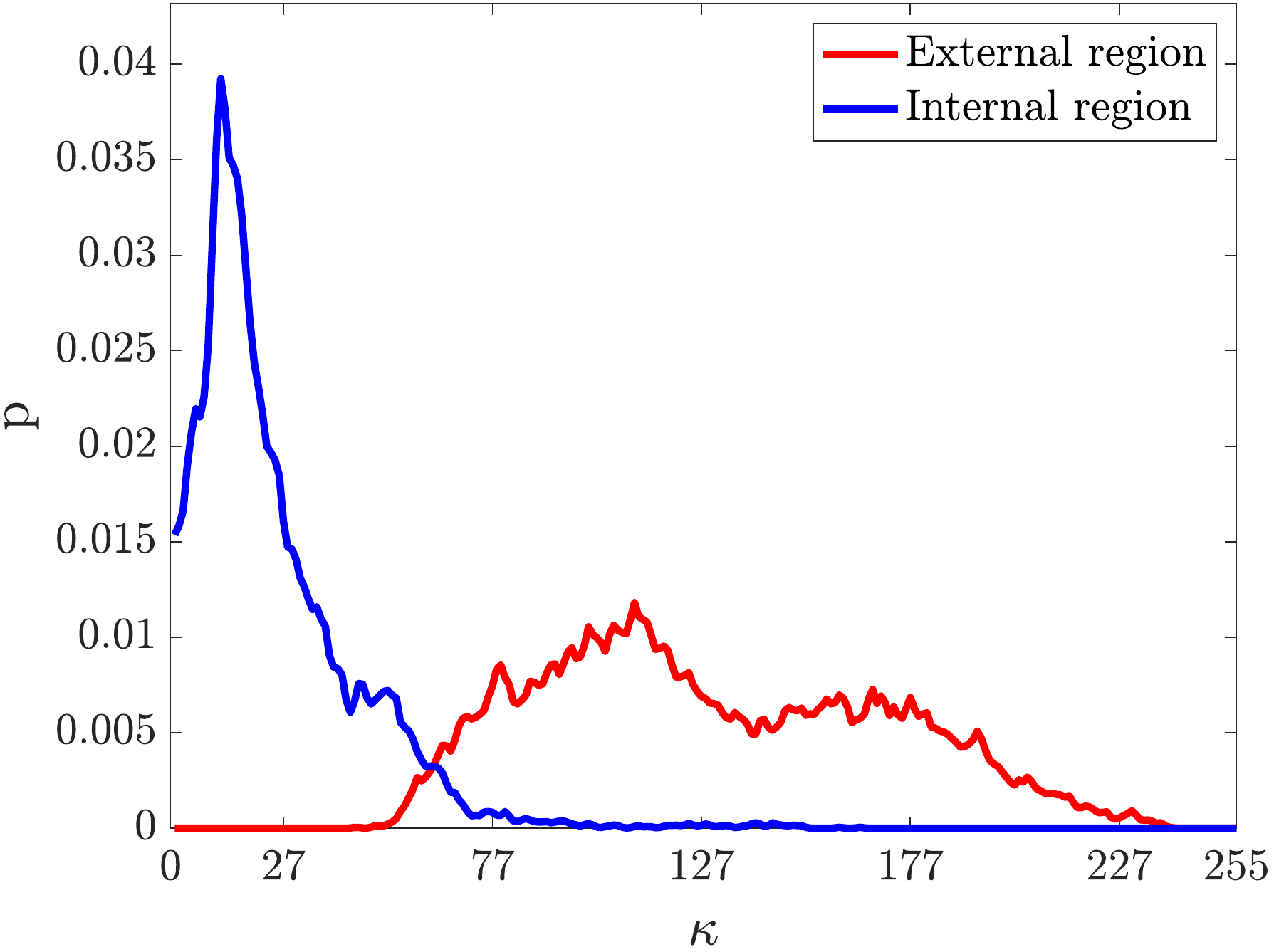}
	
\caption{Real ultrasound image - Estimated probability density function of the pixel intensity $\kappa$.}
\label{fig:ultrasoundProb}
\end{SCfigure}

The employed Bayesian segmentation model is able to extract the spatial information from the image under analysis, despite its complexity and the presence of speckle noise. This is confirmed by the PDF of the pixel intensity $\kappa$ in the internal and external regions, obtained via the nonparametric kernel density estimate (see figure~\ref{fig:ultrasoundProb}). Indeed, the model clearly identifies two regions in the image, each associated with a different interval of values of the pixel intensity.
The numerical artifacts yielded by the split-adapt Bregman segmentation can be ascribed to the high level of noise in the original image. A possible remedy to this issue, beyond the scope of the present manuscript, is the coupling of the discussed segmentation strategy with a preliminary denoising or filtering step in order to improve the quality of the input data for the split-adapt Bregman method.

\section{Concluding remarks}
\label{sc:Conclusion}

In this work, the split Bregman algorithm for region-based image segmentation was enriched by means of a mesh adaptation procedure driven by an anisotropic recovery-based error estimate. The methodology employs a piecewise constant definition of the level-set function to describe the interface between the background and the foreground of the image, leading to large values of the gradient of such a function across the contour of the region to be segmented. In addition, the image spatial information is exploited during the segmentation via the construction of a Bayesian model which estimates the probability distributions of the pixel intensity in each region.

The proposed split-adapt Bregman algorithm solves the first step of the split Bregman method by means of the finite element method. Then, an anisotropic recovery-based error estimate is computed to construct a new metric (i.e., to generate a new mesh) in order to capture the directional features of the gradient of the level-set function. The resulting anisotropically adapted mesh provides an accurate  description of the interface between the background and the foreground of the image, with a reduced number of degrees of freedom. Finally, the optimisation algorithm is resumed and the \emph{shrinkage} step and the Bregman update are performed on the adapted mesh.

\hl{The first novelty of the proposed approach is the capability to handle image inhomogeneities by exploiting the spatial information encapsulated in the PDFs of the Bayesian model. This ensures an accurate segmentation of real images featuring complex spatial patterns that traditional region-based or edge-based segmentation strategies fail to process.  In particular,  the Bayesian model is especially suited to integrate available information on the image into the model through the prior term in equation~\eqref{eq:MAP}, thus improving the segmentation performance.  Of course,  tailored energy functionals could be devised to specifically target certain classes of images, thus outperforming the presented segmentation model. This is, e.g., the case of the synthetic images with homogeneous regions tested in section~\ref{sc:Validation}, for which the RSFE model displays slightly more accurate results than the Bayesian model. Nonetheless,  the Bayesian model provides a more general approach, robust to different sources of noise.

The second novelty of this work is the integration of a mesh adaptation step in the optimisation procedure to enhance the accuracy of the segmentation. Indeed, this ensures a smooth description of the interface without jagged details.  The proper alignment, sizing and shaping of the anisotropically adapted mesh elements also guarantee that the increased precision of the contour is achieved with a reduction of the degrees of freedom up to $94 \%$, making the overall split-adapt Bregman method very competitive from a computational viewpoint.  Moreover,  the split-adapt Bregman procedure is not limited to the discussed Bayesian segmentation model, as testified by its successful application to the RSFE functional in section~\ref{sc:Validation}. As a matter of fact, the proposed optimisation-and-adaptation strategy is suitable to treat any energy functional that can be written according to the form~\eqref{eq:minL1phi}, extending the benefits of the split-adapt Bregman method to other segmentation models.
}

The resulting split-adapt Bregman method was successfully tested on synthetic and real images featuring inhomogeneous spatial patterns. In addition, the algorithm proved to be robust even in the presence of different types of noise, including Gaussian, salt and pepper and speckle noise, confirming its superiority with respect to the standard split Bregman approach. Finally, the split-adapt Bregman method was applied to an extremely challenging medical image obtained from the ultrasound of a gallbladder with remarkably good results.

\section*{Acknowledgements}
This work was developed during a research stay of M.G. at Politecnico di Milano, with the financial support of the SIMAI-ACRI Young Investigators Training Programme (YITP).
M.G. also acknowledges the support of the Spanish Ministry of Science and Innovation and the Spanish State Research Agency MCIN/AEI/10.13039/501100011033 (Grants No. PID2020-113463RB-C33 and CEX2018-000797-S) and the Generalitat de Catalunya through the Serra H\'unter Programme.
S.P. gratefully acknowledges the financial support of INdAM, Italy - GNCS Project 2020.

\bibliographystyle{elsarticle-num}
\bibliography{Ref-Img}

\appendix

\section{Region scalable fitting energy approach}
\label{sc:appRSFE}

Starting from the work~\cite{ChanVese-01} by T.F. Chan and L.A. Vese on region-based models, the RSFE method~\cite{Li-LKGD-08} improves the description of spatially inhomogeneous regions by exploiting the local information in the neighbourhood of each pixel. This is achieved by introducing the local weighted averages $f_i, \ i {=} \mathrm{I},\mathrm{E}$ of the pixel intensity given by
\begin{equation}
\label{eq:fIntExt}
f_i(\bx) :=  \frac{\displaystyle\int_{\Omega}{\Ksig(\bx-\bz) U(\bz) \chi_i(\phi(\bz)) d\bz}}{\displaystyle\int_{\Omega}{\Ksig(\bx-\bz) \chi_i(\phi(\bz)) d\bz}} ,
\end{equation}
where $\chiI(\phi) {=} \Heps(\phi)$ and $\chiE(\phi) {=} 1 {-} \Heps(\phi)$ denote the characteristic functions associated with $\OmegaI$ and $\OmegaE$, respectively, while $\Ksig$ is a Gaussian kernel with standard deviation $\sigma$.
Notice that the classical averages of the pixel intensity characterising the Chan-Vese model can be retrieved by neglecting the Gaussian kernel $\Ksig$ in equation~\eqref{eq:fIntExt}.

The discrepancy between the pixel intensity of the image $U$ and the local weighted averages $\fI$ and $\fE$ in $\OmegaI$ and $\OmegaE$, respectively are defined by
\begin{equation}
\label{eq:eIntExt}
\begin{aligned}
\eI(\bx) &:=  \int_{\Omega}{\Ksig(\bx-\bz) (U(\bz)-\fI(\bx))^2 \chiI(\phi(\bz)) d\bz} \\
\eE(\bx) &:=  \int_{\Omega}{\Ksig(\bx-\bz) (U(\bz)-\fE(\bx))^2 \chiE(\phi(\bz)) d\bz} ,
\end{aligned}
\end{equation}
where the convolution product with the Gaussian kernel is employed.

Following~\cite{Li-LKGD-08}, the energy functional for the RSFE model thus coincides with
\begin{equation}
\label{eq:funcRSFE}
\Frsf(\phi,\eI,\eE)  := \int_{\Omega}{ [ \muI \eI(\bx) + \muE \eE(\bx) ] d\bx} 
+ \nu \int_{\Omega}{g(U(\bx)) | \nabla \phi(\bx) | d\bx} ,
\end{equation}
where $\muI$, $\muE$ and $\nu$ are positive constants to be properly tuned. The first term is responsible for the separation of $\Omega$ into two regions by exploiting the discrepancies $\eI$ and $\eE$ defined in equation~\eqref{eq:eIntExt}, whereas the second term introduces a regularisation of the segmented contour via the weighted TV norm as in section~\ref{sc:Bayes}.

\subsection*{Split Bregman algorithm for region-scalable fitting energy}

Following the rationale introduced in section~\ref{sc:SplitBreg}, the minimisation of the RSFE functional~\eqref{eq:funcRSFE} can be rewritten in the form of  the optimisation problem~\eqref{eq:minL1phi}, with the functional $\G(\phi(\bx))$ defined as
\begin{equation}
\label{eq:funcGrsf}
\G(\phi(\bx)) = \Grsf(\phi(\bx)) := \int_{\Omega}{ (\muI \eI(\bx) + \muE \eE(\bx)) d\bx} .
\end{equation}
Hence, steps A, B and C of the split Bregman strategy in equation~\eqref{eq:splitBregman} can be similarly applied to the minimisation of the RSFE functional. Only minor modifications are required to adapt the framework presented in section~\ref{sc:SplitBreg} to the new energy functional in equation~\eqref{eq:funcGrsf}. More  precisely, the evolution equation~\eqref{eq:evolutionPhi} yielded by step A of the split Bregman strategy requires the definition of a new source term
\begin{equation}
\label{eq:sourceTermRSF}
s(\bx) = \sRSF(\bx) := \muI \eI(\bx) + \muE \eE(\bx) ,
\end{equation}
which is thus evaluated using the last computed level-set function $\phi^k$ as in section~\ref{sc:SplitBreg}.
The remaining steps of the split Bregman algorithm are not affected by the change in the definition of the energy functional from $\Fbay$ to $\Frsf$. Hence, equations~\eqref{eq:shrinkD} and~\eqref{eq:splitBregman-C} are valid also for the minimisation of the RSFE functional and the overall procedure described in algorithm~\ref{alg:splitBregman} stands. Finally, the optimisation procedure stops when the relative variation of the level-set function from iteration $k$ to iteration $k {+} 1$ is below a user-defined tolerance $\eta^\star$ (Algorithm~\ref{alg:splitBregman} - Step 3).

\section{Selection of model and method parameters}
\label{sc:appParam}

\hl{
In table~\ref{tab:params}, the values of the parameters involved in the definition of the energy functionals and in the construction of the split Bregman and split-adapt Bregman methods for the simulations in section~\ref{sc:Simulations} are provided. } For a detailed discussion on the sensitivity of the split Bregman method to the choice of its defining parameters, interested readers are referred  to~\cite{Li-LKGD-08,Osher-YLKO-10,Brox-BC-09}.
\begin{table}[!htb]
\begin{center}
\hl{
\begin{tabular}{ | C{3cm} | C{3cm} | C{3cm} | }
\hline
\multicolumn{3}{|c|}{Bayesian model accounting for spatial information}                                  
\\ \hline
$\tau$      & $\mu$    & $\zeta$
\\ \hline
$1$           &  $1$        & $10^{-8}$
\\ \hline
\multicolumn{3}{|c|}{Region-scalable fitting energy}                                                     
\\ \hline
$\sigma$  & $\mu$     & $\mu_i, \ i {=} \mathrm{I},\mathrm{E}$
\\ \hline
$8$           & $10^{-3}$ & $10^{-5}$
\\ \hline
\multicolumn{3}{|c|}{Split Bregman algorithm}                                                               
\\ \hline
$\nu$       & $\beta$   & $\varepsilon$
\\ \hline
$1$          & $100$     & $10^{-2}$
\\ \hline
\multicolumn{3}{|c|}{Split-adapt Bregman algorithm}                                                               
\\ \hline
$\tau^\star$   &   $\nRec$   &   $\omega$
\\ \hline
$0.5$       &       $3$      &     $0.9$
\\ \hline
\end{tabular}
}
 \end{center}
  \caption{\hl{Parameters used in the definition of the Bayesian and of RSFE functionals, and in the setup of the split Bregman and of the split-adapt Bregman algorithms.}}
    \label{tab:params}
\end{table}

\end{document}